\documentclass[11pt]{amsart}

\usepackage{amsmath,yhmath}
\usepackage{amssymb}
\usepackage{amsfonts}
\usepackage{mathrsfs} 
\usepackage[mathscr]{eucal} 
\usepackage{bbm} 
\usepackage{mathtools}
\usepackage[all]{xy}
\usepackage{tikz}
\usepackage{tikz-cd}
\usepackage[inline]{enumitem}
\allowdisplaybreaks
\usepackage{stmaryrd} 

\numberwithin{equation}{section}

\usepackage{caption}
\usepackage{geometry}
\usepackage{parskip}
\usepackage{setspace}
\usepackage[colorlinks, linkcolor=black,citecolor=blue, urlcolor=cyan, unicode, bookmarksnumbered]{hyperref}

\usepackage{aliascnt}

\newtheoremstyle{note}
{\topsep}
{\topsep}
{}
{0pt}
{\bfseries}
{.}
{0.5em }
{}
\theoremstyle{plain}

\newtheorem{Thm}{Theorem}[section]

\newtheorem*{Thm*}{Theorem}

\newaliascnt{prop}{Thm}
\newtheorem{Prop}[prop]{Proposition}
\aliascntresetthe{prop}

\newaliascnt{lemma}{Thm}
\newtheorem{Lemma}[lemma]{Lemma}
\aliascntresetthe{lemma}

\newaliascnt{coro}{Thm}
\newtheorem{Coro}[coro]{Corollary}
\aliascntresetthe{coro}

\newaliascnt{conjecture}{Thm}

\aliascntresetthe{conjecture}

\newtheorem{Lemma*}{Lemma}
\newtheorem*{Prop*}{Proposition}
\newtheorem*{Coro*}{Corollary}

\theoremstyle{definition}

\newaliascnt{def}{Thm}
\newtheorem{Def}[def]{Definition}
\aliascntresetthe{def}

\newtheorem*{Def*}{Definition}

\newaliascnt{Eg}{Thm}
\newtheorem{eg}[Eg]{Example}
\aliascntresetthe{Eg}

\newtheorem*{eg*}{Example}

\theoremstyle{remark}

\newaliascnt{rmk}{Thm}
\newtheorem{RMK}[rmk]{Remark}
\aliascntresetthe{rmk}

\newtheorem*{RMK*}{Remark}

\newcommand{\bC} { {\mathbb{C}}}

\newcommand{\bZ} { {\mathbb{Z}}}
\newcommand{\bR} { {\mathbb{R}}}

\newcommand{\bN} { {\mathbb{N}}}
\newcommand{\cA} { {\mathcal{A}}}
\newcommand{\cB} { {\mathcal{B}}}
\newcommand{\cC} { {\mathcal{C}}}
\newcommand{\cD} { {\mathcal{D}}}
\newcommand{\cE} { {\mathcal{E}}}

\newcommand{\cO} { {\mathcal{O}}}
\newcommand{\cV} { {\mathcal{V}}}

\newcommand{\sT} { {\mathscr{T}}}

\newcommand{\HOM}{{ \textnormal{Hom}  }}

\newcommand{\pt}{{\operatorname{pt}}}

\newcommand{\id}{{\operatorname{id}}}
\newcommand{\Id}{{\operatorname{Id}}}
\newcommand{\Mod}{{\operatorname{Mod}}}
\newcommand{\fpw}[1]{{\llbracket #1\rrbracket}}
\newcommand{\Sh} { {\operatorname{Sh}}}

\newcommand{\colim}{\operatornamewithlimits{colim}}
\newcommand\inEnd{\underline{\operatorname{End}}}
\newcommand\End{\operatorname{End}}
\newcommand\Fun{\operatorname{Fun}}
\newcommand\HH{\operatorname{HH}}
\newcommand{\Thom}{\operatorname{Th}}
\newcommand{\Thompr}{\mathcal{T}\mathit{h}}
\newcommand\Tr{\operatorname{Tr}}

\newcommand\Aut{\operatorname{Aut}}

\newcommand\CAlg{ {\operatorname{CAlg}} }

\newcommand\Hom{\operatorname{Hom}}

\newcommand{\Map}{{\operatorname{Map}}}

\newcommand\Pic{ \operatorname{Pic}  }
\newcommand\PrL{  {  \operatorname{Pr}^{\operatorname{L} }  } }

\newcommand\PrLst{  {  \operatorname{Pr}^{\operatorname{L} }_{st} } }

\newcommand\fib{\operatorname{fib}}

\newcommand\Sp{ {\operatorname{Sp} } }
\newcommand\An{ {\operatorname{An} } }
\newcommand\Vect{\operatorname{Vect}}

\newcommand\Loc{ {\operatorname{Loc}} }

\newcommand\sHom{\operatorname{\mathscr{H}om}}

\newcommand\Coh{\operatorname{Coh}}
\newcommand\IndCoh{\operatorname{IndCoh}}
\newcommand\Perf{\operatorname{Perf}}
\newcommand\QCoh{\operatorname{QCoh}}
\newcommand\Spec{\operatorname{Spec}}

\newcommand\HC{ {\operatorname{HC} } }

\newcommand\PrLV[1][\sT]{  {  \operatorname{Pr}^{\operatorname{L} }_{\sT,st} } }
\newcommand\Perv{\operatorname{\Perv}}

\newcommand\nT{ {\textnormal{T}} }

\newcommand{\larbR} { \overleftarrow{\mathbb{R}}}

\makeatletter
\newcommand\boxstar{\mathrel{\boxcls@{\star}}}
\newcommand{\boxcls@}[1]{%
  \vphantom{\Box}
  \mathpalette\boxcls@@{#1}%
}
\newcommand{\boxcls@@}[2]{%
  \ooalign{$\m@th#1\Box$\cr
  \hidewidth\boxcls@fix{#1}\hbox{$\m@th#1#2\mkern 2.4mu$}\cr}}

\newcommand\boxcls@fix[1]{
  \ifx#1\displaystyle
    \raise.225ex
  \else
    \ifx#1\textstyle
      \raise.225ex
    \else
      \ifx#1\scriptstyle
        \raise.180ex
      \else
        \raise.150ex
      \fi
    \fi
  \fi
}
\makeatother

\title{Cyclic homology of Sheaf Categories}
\author{Wenyuan Li}
\address{Wenyuan Li,\newline School of Mathematical Sciences, Zhejiang University, 866 Yuhangtang Rd, Hangzhou, 310058, P. R. China.}
\email{wenyuan.li@zju.edu.cn}

\author{Bingyu Zhang}
\address{Bingyu Zhang,\newline Department of Mathematics, Kyiv School of Economics, 3 Mykoly Shpaka St, Kyiv, 02000, Ukraine.}
\email{bzhang@kse.org.ua}

\curraddr{School of Mathematics, Institute for Advanced Study, 1 Einstein Dr, Princeton, NJ 08540, USA.}
\email{byzhang@ias.edu}

\date{\today}

\begin{document}

\begin{abstract}   In this note, we explain the circle action on the Hochschild homology of categories of sheaves and their admissible subcategories in terms of the six-functor formalism using factorization homology, including (quasi-)coherent sheaves on derived stacks, sheaves on manifolds with prescribed microsupports, and Tamarkin categories with prescribed microsupports. We also explain the Gysin sequence for $C_\ell$ and $S^1$-(co)invariants and the K\"unneth sequence for $S^1$-(co)invariants of the Hochschild homology relative to any complex oriented symmetric monoidal category.
\end{abstract}
\maketitle

\section{Introduction}
Hochschild homologies and cyclic homologies, introduced by Connes \cite{Connes1983}, Tsygan \cite{Tsygan-homology-Lie-algebra}, Loday and Quillen \cite{Loday-Quillen-cyclic}, are important invariants of (dg-)algebras and dg-categories \cite{McCarthy,Keller-exact} and can be generalized to ring spectra and stable categories \cite{Bokstedt-Hsiang-Madsen,Nikolaus-Scholze-Ontopologicalcyclichomology,AMGR-FH-enriched}. For categories of sheaves, the Hochschild and cyclic homologies often recover important invariants of the underlying spaces. For instance, for coherent sheaves, the Hochschild homology is isomorphic to differential forms on the (smooth) scheme by the Hochschild--Kostant--Rosenberg theorem \cite{FeigenTsygan,LodayCyclic,Schuhmacher,Yekutieli,Ben-Zvi-Nadler-looplanglands,toen_vezzosi_2011}.

In general, the categories of sheaves one studies may fail to be compactly generated (for instance, sheaves on locally compact Hausdorff spaces are rarely compactly generated \cite{Neeman,Harr}; even quasi or ind-coherent sheaves on derived stacks may fail to be compactly generated \cite{Ben-Zvi-Francis-Nadler,Drinfeld-Gaitsgory,stefanich2025dualizability_deriveds_stack}). Therefore, one needs to extend the definition of Hochschild and cyclic homology to dualizable stable categories \cite{Hoyois-Scherotzke-Sibilla,Ayala-Francis-FH-dim1,Efimov-K-theory}.

In the present note, we will record the formulas for cyclic homologies of certain categories of sheaves, using the general results on cyclic homologies of dualizable stable categories by factorization homology \cite{Ayala-Francis-FH,AMGR-FH-enriched,Ayala-Francis-FH-dim1}. We will also record the Gysin sequence and K\"unneth sequence for cyclic homologies relative to complex oriented coefficients. While we expect most results are known to experts, we think this could serve as a reference for future work (one can find relevant results in \cite{toen_vezzosi_15,Ben-Zvi-Nadler-trace,Hoyois-Scherotzke-Sibilla} which obtain similar constructions directly using \cite{Lurie-cobordism}, but it seems not immediately clear how to formulate the circle action combinatorially as a cyclic object).

\subsection{Topological sheaves}
The microlocal theory of sheaves, introduced and developed by Kashiwara and Schapira \cite{KS90}, has proven to be an important tool for studying symplectic and contact geometry \cite{nadler2009constructible,tamarkin2013}. They are closely related to the Fukaya categories \cite{nadler2009constructible,GPS3}, along with many other symplectic and contact invariants defined using pseudo-holomorphic curves.

We will assume that $\cE$ is a compactly generated symmetric monoidal stable category, for example, modules over a commutative ring or a commutative ring spectrum. For a paracompact topological space $X$, we write $\Sh(X) \coloneqq \Sh(X; \cE)$ for the category of sheaves valued in $\cE$. This is a symmetric monoidal presentable stable category.

For a manifold $M$ and a sheaf $F \in \Sh(M)$, Kashiwara and Schapira introduced the notion of singular support $ss(F) \subset T^*M$ which is a conic closed subset that measures the failure of propagations of sections along different codirections \cite{KS90}. We write $\Sh_Z(M) \coloneqq \Sh_Z(M; \cE)$ for the category of sheaves such that $ss(F) \subset Z$. 

Recently, people have written down a formula for the Hochschild homology with circle action of $\Sh_Z(M)$ when $Z \subset T^*M$ is an isotropic subset \cite{Kuo-Li-CY}. We will explain that the formula holds for a general closed subset $Z \subset T^*M$.

For a general closed subset $Z \subset T^*M$, this is a large but dualizable category, and the Hochschild homology can be defined via the trace of the identity functor \cite{Hoyois-Scherotzke-Sibilla}. Using the results in factorization homology \cite{Ayala-Francis-FH-dim1}, one can define the $S^1$-action on the Hochschild homology in an explicit way. 

It turns out that $\Sh_Z(M) \subset \Sh(M)$ is an admissible subcategory (see \autoref{sec:cat-of-sheaf}), and the projector is realized by convolution of sheaf kernels:
\[K_Z\circ -: \Sh(M) \twoheadrightarrow \Sh_Z(M).\]

We fix some notation. Let $X$ be the underlying space, $\pi_{n} : X^n  \rightarrow \pt$ be the projection, and $\widetilde{\Delta}_n$ be the twisted diagonal map:
\begin{align*}
    \widetilde{\Delta}_n&: X^n \rightarrow X^{2n}, \qquad (x_1,\cdots,x_n)\mapsto (x_n,x_1,x_1,\cdots,x_{n-1},x_{n-1},x_n).
\end{align*}

\begin{Thm}[\autoref{cor:trace-sheaf}]
    Let $Z \subset T^*M$ be a closed conic subset. The Hochschild homology with its $S^1$-action of sheaves microsupported in $Z$ is isomorphic to the cyclic object $\HH_*(\Sh_Z(M)): \boldsymbol{\Lambda} \to \cE$ defined by:
    \[\HH_*(\Sh_Z(M), n) = \pi_{n!}\widetilde{\Delta}_{n}^*(K_Z \boxtimes \dots \boxtimes K_Z).\]
    It is functorial under inclusions $j: Z_1 \hookrightarrow Z_2$ and invariant under conic Hamiltonian diffeomorphisms $\varphi: Z_1 \xrightarrow{\sim} Z_2$.
\end{Thm}

Using the framework of the $S^1$-action on Hochschild homology, we also define the $C_\ell$-action on Hochschild homology, which, for compactly generated categories, agrees with Chen \cite{Chen1} using the $C_\ell$-cyclic bar construction.

\begin{Thm}[\autoref{cor:trace-sheaf}]
    Let $Z \subset T^*M$ be a closed conic subset. The Hochschild homology with its $C_\ell$-action of the category of sheaves microsupported in $Z$ is isomorphic as $C_\ell$-modules to
    \[\HH_*(\Sh_Z(M)) = \pi_{\ell!}\widetilde{\Delta}_{\ell}^*(K_Z \boxtimes \dots \boxtimes K_Z),\]
    where the $C_\ell$-action is induced by cyclic permutation on the right hand side.
    It is functorial under inclusions $j: Z_1 \hookrightarrow Z_2$ and invariant under conic Hamiltonian diffeomorphisms $\varphi: Z_1 \xrightarrow{\sim} Z_2$.
\end{Thm}

Finally, when $\cE$ is a complex oriented symmetric monoidal category, we have the Gysin sequence and K\"unneth sequence for $S^1$ and $C_\ell$-(co)invariants of the Hochschild homology relative to $\cE$:

\begin{Thm}[\autoref{cor: Gysin sheaf}]
    Let $Z \subset T^*M$ be a closed conic subset. When $\cE$ is complex oriented, there is a Gysin sequence for the homotopy $S^1$ and $C_\ell$-quotients:
    \[\HH_*(\Sh_Z(M))_{hC_\ell} \to \HH_*(\Sh_Z(M))_{hS^1} \xrightarrow{[\ell](u)} \HH_*(\Sh_Z(M))_{hS^1}[2]\]
    and a Gysin sequence for the homotopy $S^1$ and $C_\ell$-fixed points:
    \[\HH_*(\Sh_Z(M))^{hS^1}[-2] \xrightarrow{[\ell](u)} \HH_*(\Sh_Z(M))^{hS^1} \to \HH_*(\Sh_Z(M))^{hC_\ell}\]
    where $u$ is the first Chern class in $H^2(BS^1; E)$ and $[\ell](u)$ is the $\ell$-series of the formal group law.
\end{Thm}

\begin{Thm}[\autoref{cor: Künneth sheaf}]
    Let $Z \subset T^*M$, $Y \subset T^*N$ be closed conic subsets. When $\cE$ is complex oriented, there is a K\"unneth isomorphism for the homotopy $S^1$-quotients:
    \[\HH_*(\Sh_{Z\times Y}(M \times N))_{hS^1} \simeq \HH_*(\Sh_Z(M))_{hS^1} \otimes _{E^{hS^1}} \HH_*(\Sh_Y(N))_{hS^1}[1],\]
    and a K\"unneth isomorphism for the homotopy $S^1$-fixed points within derived complete $\cE$-modules:
    \[\HH_*(\Sh_{Z\times Y}(M \times N))^{hS^1} \simeq \HH_*(\Sh_Z(M))^{hS^1} \widehat\otimes _{E^{hS^1}} \HH_*(\Sh_Y(N))^{hS^1}.\]
\end{Thm}

\subsection{Tamarkin categories}
To study quantitative symplectic geometry of $T^*X$ using sheaf theory, Tamarkin defined a category $\sT(T^*X)$ as a Verdier quotient of the category of sheaves $\Sh(X\times \bR)$ by a microlocal condition, where the extra $\bR$-factor records the action filtration \cite{tamarkin2013}. The definition can be generalized to any open subset $U \subset T^*X$, and we will call the resulting category the Tamarkin category $\sT(U)$. Recently, this was also independently studied by \cite{Scholze-betti}.

It has been shown that the Hochschild cohomology of the Tamarkin category, or the dual of Hochschild homology of the category, is computed by the Chiu--Tamarkin complex \cite{chiu2017,zhang2020quantitative,CyclicZHANG}, and is quasi-isomorphic to the symplectic cohomology for nice symplectic domains \cite{Hochschild-Kuo-Shende-Zhang}. We will explain the natural $S^1$-action on Hochschild homology, which can be used to construct symplectic capacities associated to open sets $U \subset T^*X$ \cite{zhang2020quantitative,CyclicZHANG}, analogous to the Gutt--Hutchings capacities on the $S^1$-equivariant symplectic cohomology \cite{gutthutchings2018capacities}.

Given a non-conic subset in the cotangent bundle, Tamarkin \cite{tamarkin2013} considered the lifting through the symplectic reduction diagram
\[T^*(M \times \bR) \xleftarrow{i} J^1(M) \xrightarrow{\pi} T^*M.\]
For an open subset $U \subset T^*M$, the Tamarkin category is defined by
\[\sT(U) \coloneqq \Sh(M \times \bR) / \Sh_{(\pi^{-1}(T^*M\setminus U) \times \bR_+) \cup \{\tau \leq 0\}}(M \times \bR).\]
Write $\sT = \sT(pt)$. Then the Tamarkin category $\sT(U)$ is a $\sT$-linear stable category.

This is again a large but dualizable category, and the Hochschild homology can be defined via the trace of the identity functor \cite{Hoyois-Scherotzke-Sibilla}. Using the results in factorization homology \cite{Ayala-Francis-FH-dim1}, one can define the $S^1$-action in an explicit way, which turns out to be the mixed complex defined in \cite{CyclicZHANG}.

Let $X$ denote the underlying space, $\pi_{n} : X^n  \rightarrow \pt$ be the projection, and $\widetilde{\Delta}_n$ be the twisted diagonal map:
\begin{align*}
    \widetilde{\Delta}_n&: X^n \rightarrow X^{2n}, \qquad (x_1,\cdots,x_n)\mapsto (x_n,x_1,x_1,\cdots,x_{n-1},x_{n-1},x_n).
\end{align*}

\begin{Thm}[\autoref{thm:S1-action}]
    Let $U \subset T^*M$ be an open subset. The   Hochschild homology with its $S^1$-action of the Tamarkin category is isomorphic to the $S^1$-equivariant Chiu--Tamarkin complex defined by:
    \[\HH_*^\sT(\sT(U), n) = \pi_{n!}^\sT\widetilde{\Delta}_{n,\sT}^*(P_U \boxtimes_\sT \dots \boxtimes_\sT P_U).\]
    It is functorial under inclusions $i: U_1 \hookrightarrow U_2$ and invariant under Hamiltonian homeomorphisms $\varphi: U_1 \xrightarrow{\sim} U_2$.
\end{Thm}

We also define the $C_\ell$-action on Hochschild homology of dualizable categories, and prove the $C_\ell$-Gysin sequence over any complex oriented ring spectrum, generalizing results in \cite{Chen1}.

Our second result realizes the $C_\ell$-action on Hochschild homology as a purely sheaf theoretic construction given by the $C_\ell$-equivariant Chiu--Tamarkin complex \cite{Capacities2021} motivated by \cite{chiu2017}, and compares it with the $S^1$-action on Hochschild homology over any complex oriented ring spectrum.

\begin{Thm}[\autoref{thm:Z/ell-action}]
    Let $U \subset T^*M$ be an open subset. The Hochschild homology with its $C_\ell$-action of the Tamarkin category is isomorphic to the $C_\ell$-equivariant Chiu--Tamarkin complex defined by sheaf-theoretic six-functor constructions as $C_\ell$-modules with persistence structures:
    \[\HH_*^\sT(\sT(U)) = \pi_{\ell!}^\sT\widetilde{\Delta}_{\ell,\sT}^*(P_U \boxtimes_\sT \dots \boxtimes_\sT P_U),\]
    where the $C_\ell$-action is induced by cyclic permutation on the right hand side.
\end{Thm}

when $\cE$ is a complex oriented symmetric monoidal category, we have the Gysin and K\"unneth sequences for $S^1$ and $C_\ell$-(co)invariants of the Hochschild homology relative to $\sT$ defined over $\cE$:

\begin{Thm}[\autoref{thm: gysin Tamarkin}]
    Let $U \subset T^*M$ be an open subset. When $\cE$ is complex oriented, there is a Gysin sequence between the homotopy $S^1$ and $C_\ell$ quotients
    \[\HH_*^\sT(\sT(U))_{hC_\ell} \to \HH_*^\sT(\sT(U))_{hS^1} \xrightarrow{[\ell](u)} \HH_*^\sT(\sT(U))_{hS^1}[2],\]
    and a Gysin sequence between the homotopy $S^1$ and $C_\ell$-fixed points
    \[\HH_*^\sT(\sT(U))^{hS^1}[-2] \xrightarrow{[\ell](u)} \HH_*^\sT(\sT(U))^{hS^1} \to \HH_*^\sT(\sT(U))^{hC_\ell},\]
    where $u$ is the first Chern class in $H^2(BS^1; E)$ and $[\ell](u)$ is the $\ell$-series of the formal group law.
\end{Thm}

\begin{Thm}[\autoref{thm:Künneth tamarkin}]
    Let $U \subset T^*M$, $W \subset T^*N$ be open subsets. When $\cE$ is complex oriented, there is a K\"unneth isomorphism for the homotopy $S^1$-quotients:
    \[\HH_*^\sT(\sT(U \times W))_{hS^1} \simeq \HH_*^\sT(\sT(U))_{hS^1} \otimes _{(E_{[0,\infty)})^{hS^1}} \HH_*^\sT(\sT(W))_{hS^1}[1],\]
    and a K\"unneth isomorphism for homotopy $S^1$-fixed points within derived complete $\cE$-modules over $\sT$:
    \[\HH_*^\sT(\sT(U \times W))^{hS^1} \simeq \HH_*^\sT(\sT(U))^{hS^1} \widehat\otimes _{(E_{[0,\infty)})^{hS^1}} \HH_*^\sT(\sT(W))^{hS^1}.\]
\end{Thm}

We will illustrate the strengths of these algebraic results by applying them to $S^1$-equivariant symplectic capacities of open subsets $U \subset T^*M$ in the next paper \cite{Li-Zhang26-2}.

\subsection{(Quasi-)coherent sheaves}
Since the above results are formal consequences of the trace theory for dualizable categories, we can also describe the $S^1$-action on Hochschild homology of the category of topological sheaves and (quasi)coherent sheaves using six-functors. Many results can also be found in \cite{Ben-Zvi-Nadler-trace}. In particular, this generalizes Kuznetsov's description of the Hochschild homology of admissible categories of coherent sheaves of algebraic varieties \cite{KuznetsovHochschild}.

\begin{Thm}[\autoref{cor:quasicoherent}]
    Let $X$ be a perfect derived stack over a ring or a quasi-compact derived stack with affine automorphism group over a field $\Bbbk$  of characteristic $0$. The Hochschild homology with its $S^1$-action of any admissible subcategory $\cA \subset \QCoh(X)$ can be computed by the cyclic chain complex $\HH_*(\cA): \boldsymbol{\Lambda} \to D(\Bbbk)$ where
    \[\HH_*(\cA, n) = \pi_{n*}\widetilde{\Delta}_n^*(K_\cA \boxtimes \dots \boxtimes K_\cA).\]
\end{Thm}

\begin{Thm}[\autoref{cor:coherent}]
    Let $X$ be  a locally almost finite type ind-inf-scheme over a perfect field $\Bbbk$  of characteristic $0$. The Hochschild homology with its $S^1$-action of any admissible subcategory $\cA \subset \IndCoh(X)$ can be computed by the cyclic chain complex $\HH_*(\cA): \boldsymbol{\Lambda} \to D(\Bbbk)$ where
    \[\HH_*(\cA, n) = \pi_{n*}\widetilde{\Delta}_n^!(K_\cA \boxtimes \dots \boxtimes K_\cA).\]
\end{Thm}

Consider the category $\QCoh(X)$. When $X$ is a perfect stack, our description is simply an $S^1$-enhancement of the result of Ben-Zvi--Francis--Nadler \cite[Corollary 4.12]{Ben-Zvi-Francis-Nadler} and Ben-Zvi--Nadler \cite{Ben-Zvi-Nadler-trace}. Such consideration also leads to closely related interpretations of Hochschild and cyclic homologies of $\QCoh(X)$. To\"en--Vezzosi and Ben-Zvi--Nadler \cite{Ben-Zvi-Nadler-looplanglands,Ben-Zvi-Nadlerloopspaceandconnection,toen_vezzosi_2011,toen_vezzosi_15} identified the circle action on $\HH_*(\QCoh(X))$ with the circle action on the free loop space $\Gamma({L}X, \cO_{{L}X})$, and the de Rham differential on $\mathrm{Sym}_{\cO_X}(\Omega_X[1])$. When we consider only an admissible subcategory, however, our description purely in terms of six-functors seems to have an advantage.

Finally, note that we can also deduce the Gysin sequence and K\"unneth sequence for cyclic homologies in this setting. For cyclic homologies over classical rings, the Gysin sequence and K\"unneth sequence are classical (see \cite[Section 2.2 \& 4.3]{LodayCyclic}), though the $C_\ell$-Gysin sequence is relatively new (see \cite[Lecture 7]{Kaledin} and \cite[Section 5]{Chen1}).

\subsection*{Conventions}
Throughout the paper, we use the convention that all categories are $\infty$-categories. In particular, $\An$ means the $\infty$-category of anima (or homotopy types), $\Sp$ means the stable $\infty$-category of spectra, and for $E$ a commutative ring spectrum, $D(E)$ means the stable $\infty$-category of $E$-modules.  We also use the notation $\otimes$ for the monoidal structures on these categories. Throughout the paper, we use the cohomological convention, and all cohomologies $H^*(X, E)$ mean the cochain complexes or spectra in the stable $\infty$-category.

We use $C_\ell$ to denote the finite cyclic group of order $\ell$, and $\bZ/\ell\bZ$ to denote the ring of mod $\ell$ integers. We use the convention $\boldsymbol{\Delta}_{\circlearrowleft}$ for the paracyclic category, $\boldsymbol{\Lambda}$ for the cyclic category (and $\boldsymbol{\Lambda}_\ell$ for Kaledin's $\ell$-cyclic category).

\subsection*{Acknowledgment}As this work is a spiritual continuation of BZ's thesis, in which the prototype of many results was already shown, BZ would like to give Stéphane Guillermou special thanks for his encouragement. We would like to thank Chris Kuo, who is deeply involved in the discussions at the early stage of this project. We would like to thank David Ben-Zvi, John Francis, Sheel Ganatra, Ezra Getzler, Marc Hoyois, Haoyang Liu, Nick Rozenblyum and Kyler Siegel for helpful discussions. 

BZ was supported by the Novo Nordisk Foundation grant NNF20OC0066298 and VILLUM FONDEN, VILLUM Investigator grant 37814 for this work.

\subsection*{AI disclosure}The paper is split from a project initiated in the spring of 2024. We use ChatGPT 6.1 Sol to identify potential mathematical errors and polish the language. We then revise the text by ourselves.

\section{Properties of Cyclic Spectra}\label{sec:gysin}
In this section, we prove the Gysin sequence and K\"unneth formulae for general cyclic objects and $\ell$-cyclic objects defined over complex oriented symmetric monoidal categories, which will be apply to Hochschild homologies with circle action over complex oriented symmetric monoidal categories. 

\subsection{Parametrized spectra}
Consider the category of presentable categories $\PrL$ equipped with the Lurie tensor product \cite[Section 5.3.6]{HTT} \cite[Section 4.8.1]{HA}, where the category of anima $\An$ (or equivalently the category of $\infty$-groupoids) is the tensor unit of the Lurie tensor product \cite[Example 4.8.1.20]{HA}. Hence, for any $\cE \in \PrL$, we have $\cE \otimes \An\simeq \cE$, and if $\cE$ is additionally stable, we have an exact equivalence $\cE \otimes \An\simeq \cE \otimes \Sp = \cE$ induced by the smash product \cite[Corollary 1.4.2.27 \& Example 4.8.1.23]{HA}.

For a small category $X$ and a presentably symmetric monoidal category $(\cE, \otimes, E)$, we know $\cE \in \CAlg(\PrL)$, and we can set $\Fun(X,\cE)=\cE^X$. There exists a pointwise defined symmetric monoidal structure $\otimes_X$ on $\cE^X$. We have a functor
\[\An^X \times \cE \rightarrow \cE^X,\quad (F,e)\mapsto F\otimes e \coloneqq [\sigma \mapsto  F(\sigma)\otimes e]. \]
The functor preserves colimits on both factors since it is pointwise defined (as well as the colimit in functor categories), and thus descends to a functor $\An^X \otimes \cE \rightarrow \cE^X$.

\begin{Prop}[{\cite[Proposition 4.8.1.17]{HA}}]\label{prop:base-change}
Let $X$ be a small category and $\cE$ be a presentably symmetric monoidal category. Then the natural functor $\An^X \otimes \cE \to \cE^X$ defines a symmetric monoidal equivalence functorial with respect to $X$.
\end{Prop}

The straightening/unstraightening equivalence \cite[Theorem 3.2.0.1]{HTT} implies that functors from $X$ into $\An$ (resp.~$\cE$) are equivalent to left fibrations of $\An$ (resp.~$\cE$) over the simplicial set $X$:
\begin{equation}\label{eq:straighten}
\An^X \simeq \An_{/X}, \quad \cE^X \simeq \cE_{/X}.
\end{equation}

We will consider the case when $X$ is a Kan complex. Categories of parametrized spectra over a Kan complex $X$ have six functors, initiated in the works of \cite{AndoBlumbergGepner,AndoBlumbergGepnerHopkinsRezk}, which fit into a modern formalism of six-functor formalism in the sense of \cite[Definition A.5.7]{Mann} following from \cite[Proposition A.5.10]{Mann} as explained in \cite[Section 3]{HesselholtFF}. Let $f: X \to Y$ be a map between Kan complexes $X$ and $Y$. Then we have a series of adjunctions $f_! \dashv f^! = f^* \dashv f_*$ \cite[Section 5.5]{AndoBlumbergGepner} \cite[Section 3]{HesselholtFF}. They satisfy the proper base change formula. Then, consider the parametrized spectra over $\cE$, by the non-abelian base change formula \cite[Theorem 0.5]{HaineBasechange}, we have:

\begin{Prop}[{\cite[Proposition 5.12]{AndoBlumbergGepner} \cite[Theorem 0.5]{HaineBasechange}}]\label{lem:base-change}
Consider the cartesian square of Kan complexes
\[\begin{tikzcd}
    X \ar[r,"f"] \ar[d,"g"] & Y \ar[d,"h"] \\
    Z \ar[r,"i"] & W.
\end{tikzcd}\]
Then there is a canonical equivalence of functors $g_!f^* \simeq i^*h_!$ and $i^*h_* \simeq g_*f^*$ between the categories of parametrized spectra over $\cE$.
\end{Prop}

Let $X$ be a topological space. Following \cite[Section 2.6]{AndoBlumbergGepnerHopkinsRezk}, we notice that the Thom spectrum of a vector bundle $\zeta \to X$ can be defined using a parametrized spectrum over $X$ as follows. Let $\operatorname{Vect}(X)$ be the anima of vector bundles over a topological space $X$. Consider the canonical equivalence \cite[Theorem 7.1.0.1]{HTT} composed with the $J$-homomorphism
\[\operatorname{Vect}(X) \xrightarrow{\sim} \Map(X, |\Vect|) \to \Map(X, \Pic(\mathbb{S})).\]
Then we denote the image of the vector bundle $\zeta \to X$ under the composition as $\Thompr(\zeta)$. It follows from \cite[Definition 2.20]{AndoBlumbergGepnerHopkinsRezk} that the compactly supported global section of the parametrized spectrum is simply the Thom spectrum, i.e.~for the projection map $\pi: X \to pt$,
\[\Thom(\zeta) = \pi_!\Thompr(\zeta).\]

Given a fibration $f: Y \to X$, we are able to pull back parametrized spectra and therefore pullback the Thom spectrum $\Thompr(\zeta)$ to $f^*\Thompr(\zeta)$. The following result shows the relation between the two constructions:

\begin{Lemma}\label{lem:thom-pull-back}
Let $f: Y \to X$ be a map of homotopy types and $\zeta \to X$ be a vector bundle.
Then $\Thompr(f^*\zeta) = \pi_{Y!}f^*\Thompr(\zeta)$ where $\pi_Y:Y\to \pt$.
\end{Lemma}
\begin{proof}
For the vector bundle $\zeta \to X$ with corresponding classifying map $\zeta: X \to |\Vect|$, the parametrized Thom spectrum of $\zeta$ is defined by the composition
\[X \xrightarrow{\zeta} |\Vect| \xrightarrow{J} \Pic(\mathbb{S}).\]
For a map $f: Y \to X$, the pull-back vector bundle $f^*\zeta \to Y$ is by definition given by $Y \to X \to |\Vect|$. Hence, the parametrized Thom spectrum of $f^*\zeta$ is the composition functor
\[Y \xrightarrow{f} X \xrightarrow{\zeta} |\Vect| \xrightarrow{J} \Pic(\mathbb{S}).\]
Therefore, $\Thompr(f^*\zeta) = f^*\Thompr(\zeta)$.
We have $\Thompr(f^*\zeta) = \pi_{Y!}\Thompr(f^*\zeta) = \pi_{Y!}f^*\Thompr(\zeta)$.
\end{proof}

In this paper, we are particularly interested in the case that $X = BG$ for an $\infty$-group $G$ comes from a compact Lie group. In fact, we will only consider $G = C_\ell$ or $S^1$ in this article. 

The following notion of (Farrell)-Tate cohomology for topological groups is originally proposed by \cite{Klein_2001}. We follow the exposition of \cite[Section I.4]{Nikolaus-Scholze-Ontopologicalcyclichomology} and references therein.
\begin{Def}
Let $\cE$ be a presentable stable category and $M \in \cE^{BG}$. Then for $\pi: BG \to pt$, we define the homotopy orbits and fixed points as
\[M_{hG} = \pi_!M = \varinjlim\nolimits_{BG}M, \quad M^{hG} = \pi_*M = \varprojlim\nolimits_{BG}M.\]
We define the homotopy (Farrell)-Tate construction as
\[M^{tG} = \mathrm{cofib}(((\Sigma_+^\infty G)^{hG}\otimes M)_{hG} \to M^{hG}) .\]
\end{Def}

\begin{eg}If $G= C_\ell$, then $(\Sigma_+^\infty G)^{hG}=\mathbb{S}$ and the definition of $M^{tG} = \mathrm{cofib}(M_{hG} \to M^{hG}) $ recovers the usual definition of Tate cohomology. If $G=S^1$, then $(\Sigma_+^\infty S^1)^{hS^1}\simeq \Sigma \mathbb{S}$ and $M^{tS^1} = \mathrm{cofib}(\Sigma M_{hG} \to M^{hG}).$
\end{eg}

\begin{RMK}
We remark that parametrized spectra over $BG$ are different from the category of $G$-equivariant spectra used in equivariant stable homotopy theory. See for example \cite[Section II.2]{Nikolaus-Scholze-Ontopologicalcyclichomology} for their distinction and relation.
\end{RMK}

\subsection{Cyclic spectra}

Recall the combinatorial model of Connes' cyclic category $\boldsymbol{\Lambda}$, whose objects are $[n]=\{1,\dots,n\} $ for $n\in \bN_{\geq 1}$  and morphisms are generated by
\begin{align*}
    d_i: [n] \to [n+1], \;s_i: [n+1] \to [n], \;t_{n}: [n] \to [n]
\end{align*}
such that the face morphism $d_i$, the degeneration morphism $s_i$ and the cyclic rotation $t_n$ satisfy the relations in the cyclic category (with the relation $t_n^{n} = \id$) \cite[Definition 6.1.1]{LodayCyclic}\footnote{For our applications, our convention differs from the convention of \cite[Definition 6.1.1]{LodayCyclic} for the cyclic category, and then also for the simplicial category. The difference is harmless.}. Also recall the para-cyclic category $\boldsymbol{\Delta}_\circlearrowleft$, whose objects are $n \in \bN$ and morphisms are generated by
\begin{align*}
    d_i: [n] \to [n+1], \;s_i: [n+1] \to [n], \;t_{n}: [n] \to [n]
\end{align*}
such that the face morphism $d_i$, the degeneration morphism $s_i$ and the cyclic rotation $t_n$ satisfy the relations in the para-cyclic category (without the relation $t_n^n = \id$) \cite[Section 1]{Getzler-Jones} \cite[Lecture 4]{Kaledin}.

\begin{Def}
    Let $\boldsymbol{\Lambda}$ be the cyclic category. Then a cyclic object with value in $\cE$ is an object in $\Fun(\boldsymbol{\Lambda}, \cE)$. Similarly, a para-cyclic object with value in $\cE$ is an object in $\Fun(\boldsymbol{\Delta}_\circlearrowleft, \cE)$. Let $\Fun_\simeq (\boldsymbol{\Lambda}, \cE)$ and $\Fun_\simeq(\boldsymbol{\Delta}_\circlearrowleft, \cE)$ be the subcategories of functors that send all morphisms in $\boldsymbol{\Lambda}$ and $\boldsymbol{\Delta}_\circlearrowleft$ to equivalences in $\cE$, we call these objects as cartesian cyclic objects and cartesian para-cyclic objects.
    We write 
    \[()^\#: \Fun (\boldsymbol{\Lambda}, \cE) \to \Fun_\simeq (\boldsymbol{\Lambda}, \cE), \quad \text{resp.} \quad \Fun (\boldsymbol{\Delta}_\circlearrowleft, \cE) \to \Fun_\simeq(\boldsymbol{\Delta}_\circlearrowleft, \cE)\]
    for the natural localization functors.
\end{Def}

\begin{Def}Let $\cC$ be a small $\infty$-category. Let $W$ be the set of all edges in $\cC$, we call the functor $\cC\to  \cC[W^{-1}]\eqqcolon \widetilde{C}$ the $\infty$-groupoid completion.
\end{Def}
\begin{Prop}[{\cite[Proposition 6.3.1.21]{kerodon}}]\label{prop: groupoid completion}Let $\cC$ be a small $\infty$-category. We define its geometric realization as $|\cC|$ as the colimit of the constant diagram $\cC\to \An$, $c\mapsto \pt$. Then the canonical functor $\cC\to |\cC|$ is equivalent to the $\infty$-groupoid completion of $\cC$.
\end{Prop}

\begin{Prop}\label{prop:paracyclic}
    Let $\widetilde{\boldsymbol{\Lambda}}$ be the groupoid completion of $\boldsymbol{\Lambda}$ and $\widetilde{\boldsymbol{\Delta}}_\circlearrowleft$ be the groupoid completion of $\boldsymbol{\Delta}_\circlearrowleft$. Then there are $S^1$-equivariant fiber sequences of categories
    \[S^1 \to \boldsymbol{\Delta}_\circlearrowleft \to \boldsymbol{\Lambda}, \;\; S^1 \to \widetilde{\boldsymbol{\Delta}}_\circlearrowleft \to \widetilde{\boldsymbol{\Lambda}}\]
    such that $(\boldsymbol{\Delta}_\circlearrowleft)_{hS^1} = \boldsymbol{\Lambda}$ and $(\widetilde{\boldsymbol{\Delta}}_\circlearrowleft)_{hS^1} = \widetilde{\boldsymbol{\Lambda}}$. In addition, $\widetilde{\boldsymbol{\Delta}}_\circlearrowleft = pt$ and $\widetilde{\boldsymbol{\Lambda}} = BS^1$. Hence
    \[\Fun_\simeq(\boldsymbol{\Delta}_\circlearrowleft, \cE)^{hS^1} \simeq \Fun_\simeq(\boldsymbol{\Lambda}, \cE) \simeq \Fun(\widetilde{\boldsymbol{\Lambda}}, \cE) \simeq \Fun(BS^1, \cE),\]
    an object in $\cE$ with an $S^1$-action is equivalent to a cartesian cyclic object and is also equivalent to an $S^1$-invariant cartesian para-cyclic object with value in $\cE$.
\end{Prop}
\begin{proof}By \autoref{prop: groupoid completion}, the natural functor $\Lambda\to |\Lambda|$ induces
\[\Fun(|\boldsymbol{\Lambda}|,\cC) \xrightarrow{\simeq}\Fun_\simeq (\boldsymbol{\Lambda}, \cC).\] 

On the other hand, it is known that $|\boldsymbol{\Lambda}|\simeq BS^1$, for example \cite[Théorème 10]{Connes1983}.
\end{proof}

\begin{RMK}
    One may observe that the category of cartesian cyclic objects are just the category of parametrized spectra $\cE^{BS^1} = \Fun(BS^1, \cE)$ \cite{AndoBlumbergGepner,AndoBlumbergGepnerHopkinsRezk}. 
    
    We also have the straightening/unstraightening equivalence (\cite[Theorem 2.2.1.2]{HTT}) in this setting
\[\An^{BS^1} \simeq \An_{/BS^1}, \quad \cE^{BS^1} \simeq \cE_{/BS^1}.\]
Under the equivalence, a cartesian cyclic space $M \in \Fun(BS^1, \An) = \An^{BS^1}$ corresponds to the cartesian fibration $f: M_{hS^1} \to BS^1$. In particular, the point with trivial $S^1$-action represented by the parametrized space $1_{BS^1}$ corresponds to the cartesian fibration $\id: BS^1 \to BS^1$.
\end{RMK}

We can also consider $C_\ell$-actions, encoded by the category $\boldsymbol{\Delta}_{C_\ell}$, whose object is the simplex $[\ell]$ and morphism is generated by the cyclic rotation $t_\ell: [\ell] \to [\ell]$ with $t_\ell^\ell = \id$.

\begin{Def}
    A $C_\ell$-cyclic object (equivalently a cartesian $C_\ell$-cyclic object) with value in $\cE$ is an object in $\Fun(\boldsymbol{\Delta}_{C_\ell}, \cE) \simeq \Fun_\simeq(\boldsymbol{\Delta}_{C_\ell}, \cE)$.
\end{Def}

\begin{Lemma}
    There is an equivalence $\boldsymbol{\Delta}_{C_\ell} = BC_\ell$. Hence
    \[\Fun(\boldsymbol{\Delta}_{C_\ell}, \cE) \simeq \Fun_\simeq(\boldsymbol{\Delta}_{C_\ell}, \cE) = \Fun(BC_\ell, \cE),\]
    an object in $\cE$ with a $C_\ell$-action is equivalent to a $\ell$-cyclic object.
\end{Lemma}

    Consider the natural inclusion of the $C_\ell$-cyclic category to the cyclic category $\boldsymbol{\Delta}_{C_\ell} \to \boldsymbol{\Lambda}$. We have the following natural commutative diagram:

\begin{Lemma}\label{lem:finite-action-1}
    Consider the inclusion $j: \boldsymbol{\Delta}_{C_\ell} \to \boldsymbol{\Lambda}$ and the standard inclusion $j: C_\ell \to S^1$. There is a commutative diagram of categories
    \[\xymatrix{
    \boldsymbol{\Delta}_{C_\ell} \ar[r] \ar[d] & \boldsymbol{\Lambda} \ar[d] \\
    BC_\ell \ar[r] & BS^1.
    }\]
\end{Lemma}
\begin{proof}
    Let $\widetilde{\boldsymbol{\Lambda}}[\ell] \coloneqq \widetilde{\boldsymbol{\Lambda}}(-, [\ell])$ be the cyclic set represented by $[\ell]$. We know that the geometric realization of $\widetilde{\boldsymbol{\Lambda}}[\ell]$ is the automorphism group $S^1$ by  \autoref{prop:paracyclic} or \cite[Exercise 7.2.1]{LodayCyclic}. Now, it suffices to show that the geometric realization of the map $\boldsymbol{\Delta}_{C_\ell}(-, [\ell]) \to j^*\boldsymbol{\Lambda}(-, [\ell])$ is homotopic to the map on automorphism groups $C_\ell \hookrightarrow S^1$. Since the geometric realization defines a bijection on the $0$-simplices $\widetilde{\boldsymbol{\Delta}}_{C_\ell}([0], [\ell]) \to \widetilde{\boldsymbol{\Lambda}}([0], [\ell])$, it follows that it induces the inclusion $C_\ell \hookrightarrow S^1$ (see a more complicated version in \cite[Lemma 5.1]{Chen1}).
\end{proof}

Given any $S^1$-action, we can naturally obtain another $S^1$-action by considering the $\ell$-fold cover. Combinatorially, this can be described by the $\ell$-cyclic categories. Recall the $\ell$-cyclic category $\boldsymbol{\Lambda}_\ell$, whose objects are $n \in \bN$ and morphisms are generated by
\begin{align*}
    d_i: [n] \to [n+1], \;s_i: [n+1] \to [n], \;t_{n}: [n] \to [n]
\end{align*}
such that the face morphism $d_i$, the degeneration morphism $s_i$ and the cyclic rotation $t_n$ satisfy the relations in the $\ell$-cyclic category (with the relation $t_n^{nl} = \id$) \cite[Lecture 7]{Kaledin}.

\begin{Prop}\label{prop:ellcyclic}
    Let $\widetilde{\boldsymbol{\Lambda}}_\ell$ be the groupoid completion of $\boldsymbol{\Lambda}_\ell$. Then there exists a $C_\ell$-equivariant fiber sequence of categories
    \[BC_\ell \to \boldsymbol{\Lambda}_\ell \to \boldsymbol{\Lambda}\]
    such that $\boldsymbol{\Lambda}_\ell/{C_\ell} = \boldsymbol{\Lambda}$ and $\widetilde{\boldsymbol{\Lambda}}_\ell/{C_\ell} = \widetilde{\boldsymbol{\Lambda}}$. In addition, there exists an inclusion functor $\boldsymbol{\Lambda}_\ell \to \boldsymbol{\Lambda}, [n] \to [n\ell]$ that induces an equivalence $\widetilde{\boldsymbol{\Lambda}}_\ell \simeq \widetilde{\boldsymbol{\Lambda}} \simeq BS^1$.
\end{Prop}

\subsection{Complex oriented spectra}
We recall the definition of a complex oriented rigid symmetric monoidal category $\cE$. The complex orientation allows one to define the Chern class which will play an important role later for properties of cyclic objects.

Let $(\cE, \otimes, E)$ be a presentably symmetric monoidal stable category. Then we know that $\cE \in \CAlg(\PrLst)$. Let $X$ be a homotopy type or equivalently an anima. We will write the cohomology with coefficient in $\cE$ as
\[H^*(X; E) \coloneqq \Hom(\Sigma^\infty_+X,  E),\]
and the homology with coefficient in $\cE$ as
\[H_*(X; E) \coloneqq \Sigma^\infty_+X \otimes E,\]
where $\Sigma^\infty_+: \An \to \Sp$ is the infinite suspension functor.

\begin{Def}
Let $(\cE, \otimes, E)$ be a rigid symmetric monoidal stable category. Then it is complex oriented if there exists a Thom class $u \in \widetilde{H}^2(BS^1, E)$ such that the pull-back map
\[\widetilde{H}^2(BS^1, E) \to \widetilde{H}^2(S^2, E) \xrightarrow{\sim} H^0(pt, E)\]
sends $u$ to $1 \in {H}^0(pt, E)$.
\end{Def}
\begin{RMK}
One can show that complex oriented symmetric monoidal stable categories are just symmetric monoidal stable categories linear over $MU$. When $E$ is a complex oriented commutative ring spectrum, $\cE = D(E)$ is also a complex oriented symmetric monoidal category. This includes the case of complex bordism $MU$, complex $K$-theory $KU$, and all the discrete commutative rings $\Bbbk$.
\end{RMK}

It is well known that complex orientation determines a first Chern class (and in fact all Chern classes) in the cohomology theory, cf.~\cite[Proposition 6]{Lurie-chromatic}:

\begin{Lemma}\label{def:complex orientation}
Let $(\cE, \otimes, E)$ be a rigid symmetric monoidal complex oriented stable category. Then for any complex line bundle $\zeta \to X$, there exists a Thom class $u_\zeta$ such that for the Thom spectrum $\mathrm{Th}(\zeta)$,
\[u_\zeta: \Thom(\zeta) \otimes E \to \Sigma^{\infty+2}_+X \otimes E\]
is an equivalence in $\cE$.
\end{Lemma}

Let $X = BS^1$ and $\zeta_X: X \to BS^1$ be the canonical complex vector bundle defined by the identity map. Then the above isomorphism implies the Gysin sequence:

\begin{Lemma}\label{lem:gysin}
Let $\cE$ be a complex oriented symmetric monoidal stable category with a complex orientation $u$. Then there is a cofiber sequence in $\cE$
\begin{equation}\label{eq:gysin-sequence}
E \to \Sigma^\infty_+BS^1 \otimes E \xrightarrow{u} \Sigma^{\infty+2}_+BS^1 \otimes E.
\end{equation}
\end{Lemma}
\begin{proof}
Let $X = \bC P^\infty$ and $\zeta:\cO_{\bC P^\infty}(1)\rightarrow \bC P^\infty$ be the dual tautological line bundle. The underlying homotopy types realizes $X = BS^1$ and the canonical line bundle $\zeta: X \to BS^1$. The Thom spectrum is $\Thom(\cO_{\bC P^\infty}(1))=\operatorname{cofib}(\Sigma^\infty_+{pt} \to \Sigma^\infty_+\bC P^{\infty+1})$, which is the inductive limit of 
\begin{equation}\label{eq:thom}
\Thom(\cO_{\bC P^n}(1))=\operatorname{cofib}(\Sigma^\infty_+{pt} \to \Sigma^\infty_+\bC P^{n+1}).
\end{equation}
When $\cE$ is complex oriented, this implies the fiber sequence.
\end{proof}

For the following Künneth formula, we introduce the notion of formal group laws. Let $\zeta$ be the tautological line bundle over $\bC P^\infty$, $\zeta_i$ be the tautological line bundle over the $i$-th copy of $\bC P^\infty\times \bC P^\infty$, and let $\mu:\bC P^\infty\times \bC P^\infty \to BS^1$ the classifying map of  $\zeta_1\boxtimes\zeta_2$, we set
\[F(u_1,u_2)=c_1^E(\zeta_1\boxtimes\zeta_2)=\mu^*u\in H^*(\bC P^\infty\times \bC P^\infty,E)\simeq H^*(E)\fpw{u_1,u_2}.\]
Then $F$ is a formal group law in the sense of \cite{Quillen1969FormalGroupLaws}. We then define the $\ell$-series of $F$ as
\[[\ell]_E(u)=[\ell](u)\coloneqq c_1^E(\zeta^{\otimes \ell}) \in H^*(E)\fpw{u},\]
which can be given inductively using the formal group law.

\begin{Lemma}\label{lem:Künneth}
Let $\cE$ be a complex oriented symmetric monoidal stable category with a complex orientation $u$. Then there is a cofiber sequence in $\cE$
\begin{equation}
E \otimes \Sigma^\infty_+BS^1 \xrightarrow{\Delta} E \otimes \Sigma^\infty_+(BS^1 \times BS^1) \xrightarrow{u \otimes 1 - 1 \otimes u} \Sigma^2E \otimes \Sigma^\infty_+(BS^1 \times BS^1).
\end{equation}
\end{Lemma}
\begin{proof}
Let $X \times X = \bC P^\infty \times \bC P^\infty$ and $\zeta_{X \times X}: \cO_{\bC P^\infty}(1) \boxtimes \cO_{\bC P^\infty}(-1) \to \bC P^\infty \times \bC P^\infty$ be a complex line bundle, whose fiber at $(l_1,l_2)$ is isomorphic to $\Hom_\bC(l_1,l_2)$. Its classifying map is $B\chi$ for $\chi: S^1\times S^1\to S^1, (z_1,z_2)\mapsto z_1z_2^{-1}$, whose kernel is the diagonal map $\delta:S^1\to S^1\times S^1$.

Therefore, the (homotopy) fiber of $B\chi$ and the natural map $\fib(B\chi)\to BS^1 \times BS^1$ is given by the diagonal map $\Delta: BS^1 \to BS^1 \times BS^1$. 
\[\begin{tikzcd}
BS^1 \arrow[d,"\Delta"'] \ar[dr, phantom, "\square"]\arrow[r]        & ES^1 \arrow[d] \\
BS^1 \times BS^1 \arrow[r, "B\chi"] & BS^1  .        
\end{tikzcd}\]
The homotopy fiber is by definition equivalent to the sphere bundle of $\zeta_{\bC P^\infty \times \bC P^\infty}$ over $\bC P^\infty \times \bC P^\infty$.

Then the Thom spectrum is
\begin{align}\label{eq:thom-product}
    \Thom(\cO_{\bC P^\infty}(1) \boxtimes \cO_{\bC P^\infty}(-1))  
    \simeq \operatorname{cofib}(\Sigma^\infty_+\bC P^\infty \xrightarrow{\Delta} \Sigma^\infty_+(\bC P^\infty \times \bC P^\infty)). \nonumber
\end{align}
Moreover, when $u$ is the Thom class of the dual tautological line bundle $\cO_{\bC P^\infty}(1)$, the class $u \otimes 1 - 1 \otimes u$ is differs from the Thom class of the line bundle $\cO_{\bC P^\infty}(1) \boxtimes \cO_{\bC P^\infty}(-1)$, i.e. $F(u_1,[-1](u_2))$, by a unit. When $\cE$ is complex oriented, this implies the fiber sequence.
\end{proof}

\begin{RMK}
The condition that $E$ is complex oriented is also equivalent to the existence of an isomorphism 
\[E \otimes \Sigma^\infty_+BS^1 = E[u],\]
where $u$ is the Thom class of the dual tautological line bundle. Alternatively, one can also show the Gysin sequence and K\"unneth sequence using this isomorphism.
\end{RMK}

\subsection{Gysin sequences for cyclic spectra}\label{section: Gysin in general}
We prove the corresponding Gysin sequences for cyclic objects over any complex oriented symmetric monoidal category $(\cE, \otimes, E)$. Readers that are not familiar with complex oriented ring spectra may simply consider cyclic objects over $\bZ$ throughout this section.

Let $M$ be a cartesian cyclic object over $E$. By \autoref{prop:paracyclic}, we know that $M$ is a locally constant presheaf, or equivalently, a parametrized spectrum over $BS^1$. Write $j: pt \hookrightarrow BS^1$ and $\pi: BS^1 \to pt$. Following \cite[Example 6.1.6.2]{HA} and \cite[Example I.1.12]{Nikolaus-Scholze-Ontopologicalcyclichomology}, we write $j^*$ for the forgetful functor, $\pi_!$ for the homotopy orbit functor, and $\pi_*$ for the homotopy fixed point functor:
\begin{align*}
    M_{hS^1} &= \pi_!M = (\Sigma^\infty_+ES^1 \otimes M)_{S^1} =\varinjlim\nolimits_{BS^1}M, \\
    M^{hS^1} &= \pi_*M = \Hom(\Sigma^\infty_+ES^1, M)^{S^1}=\varprojlim\nolimits_{BS^1}M.
\end{align*}

Consider the trivial $S^1$-action on $E$ so that $E_{hS^1} = (\Sigma^\infty_+ES^1 \otimes E)_{S^1} = \Sigma^\infty_+BS^1 \otimes E$ and $E^{hS^1} = \Hom(\Sigma^\infty_+ES^1, E)^{S^1} = \Hom(\Sigma^\infty_+BS^1, E)$. With the notation, we have the following reformulation of Gysin sequences in \autoref{lem:gysin}
\[E \to E_{hS^1} \xrightarrow{u} E_{hS^1}[2], \;\; E^{hS^1}[-2] \xrightarrow{u} E^{hS^1} \to E.\]
We can generalize it to parametrized spectra and obtain the following Gysin sequences:

\begin{Prop}\label{prop:Gysin}
    Let $M$ be a cyclic object over a complex oriented symmetric monoidal category $(\cE, \otimes, E)$. Then there is a fiber sequence for the homotopy orbits
    \[j^*M \to M_{hS^1} \xrightarrow{u} M_{hS^1}[2].\]
    Similarly, there is a fiber sequence for homotopy fixed points
    \[M^{hS^1}[-2] \xrightarrow{u} M^{hS^1} \to j^*M.\]
\end{Prop}
\begin{proof}
    By \autoref{lem:thom-pull-back} and \autoref{eq:thom}, the Thom spectrum of the pullback bundle of $\cO_{\bC P^\infty}(1) \to \bC P^\infty$ is
    \[\Thom(f^*\cO_{\bC P^\infty}(1)) = \operatorname{cofib}(\Sigma^\infty_+j^*M \to \Sigma^\infty_+M_{hS^1}).\]
    Therefore, since $E$ is the unit of the complex oriented category $\cE$, the Thom isomorphism implies the Gysin sequence
    \[E \otimes \Sigma^\infty_+j^*M \to E \otimes \Sigma^\infty_+M_{hS^1} \xrightarrow{u} E \otimes \Sigma^\infty_+M_{hS^1}[2].\]
    The result follows from the equivalence $\cE^{BS^1} \simeq \cE \otimes \An^{BS^1}$.
\end{proof}

\begin{RMK}    One can also show the Gysin sequence using the isomorphism 
    \[H_*(BS^1, E) = E[u^{-1}]\]
    Consider the fibration $M_{hS^1} \rightarrow BS^1$. The Atiyah–Hirzebruch spectral sequence implies $H_*(j^*M) \otimes H_*E[u^{-1}] \Rightarrow H_*(M_{hS^1})$. Consider the fiber sequence 
    \begin{equation}E \to E[u^{-1}] \xrightarrow{u} E[u^{-1}][2].\end{equation}
    Using functoriality of the spectral sequences, we then obtain the Gysin sequence.
\end{RMK}

\begin{RMK}
    We can also obtain the above result using explicit computations for discrete rings. Let $M$ be a cyclic object over $\Bbbk$. Consider the cellular structure on $S^1$ with one $0$-cell and one $1$-cell, so that $S_1 = C_{-*}^{cell}(S^1)$ is represented by the algebra of dual numbers $\Bbbk[\varepsilon]$. Then there is an equivalence
    \[\Fun(BS^1, D(\Bbbk)) \xrightarrow{\sim} D(\Bbbk[\varepsilon]),\]
    such that
    \[M_{hS^1} = \pi_!M = \Bbbk \otimes_{\Bbbk[\varepsilon]} M,\qquad M^{hS^1} = \pi_*M = \Hom_{\Bbbk[\varepsilon]} ( \Bbbk,M).\]
    Consider the $\Bbbk[\varepsilon]$-module $\Bbbk$ with the resolution
    \[\dots \rightarrow \Bbbk[\varepsilon] \xrightarrow{\varepsilon} \Bbbk[\varepsilon] \xrightarrow{\varepsilon} \Bbbk[\varepsilon] \rightarrow 0.\]
    This implies that $\Bbbk \otimes_{\Bbbk[\varepsilon]} \Bbbk = \Bbbk[u^{-1}]$, $\Hom_{\Bbbk[\varepsilon]}(\Bbbk, \Bbbk) = \Bbbk\fpw{u}, \, |u| = 2$. There is a Gysin map induced by $u$:
    \[\begin{tikzcd}
{\Bbbk[\varepsilon]} \arrow[r, "{\varepsilon}"] \arrow[d] & {\Bbbk[\varepsilon]} \arrow[r, "{\varepsilon}"] \arrow[d] & {\Bbbk[\varepsilon]} \arrow[r, "{\varepsilon}"] \arrow[d, "0"] & {\Bbbk[\varepsilon]} \arrow[r, "{0}"] \arrow[d] & 0 \arrow[d] \\
{\Bbbk[\varepsilon]} \arrow[r, "{\varepsilon}"]           & {\Bbbk[\varepsilon]} \arrow[r, "{0}"]                             & 0 \arrow[r]                                                          & 0 \arrow[r]                                                       & 0  .     
\end{tikzcd}\]
    Then, by tensoring with $M$ over $\widetilde{\boldsymbol{\Lambda}}$, we get a resolution for $\Bbbk \otimes_{\Bbbk[\varepsilon]} M$, respectively $\Hom_{\Bbbk[\varepsilon]}(\Bbbk, M)$ and obtain the Gysin sequence.
\end{RMK}

Then we consider the $C_\ell$-Gysin sequence that relates the $S^1$-orbits and the $C_\ell$-orbits, or the corresponding fixed points. We recall the notions of the $C_\ell$-cyclic category $\boldsymbol{\Delta}_{C_\ell}$ which encodes $C_\ell$-actions, and the $\ell$-cyclic category $\boldsymbol{\Lambda}_\ell$ of Kaledin which encodes the $\ell$-fold cover of $S^1$-actions.
\begin{Prop}\label{prop:Gysin-finite-group}
    Let $M$ be a cyclic object over a complex oriented symmetric monoidal category $(\cE, \otimes, E)$. Let $j_\ell: \boldsymbol{\Delta}_{C_\ell} \to \boldsymbol{\Lambda}$ be the inclusion and $[\ell](u)$ be the   $\ell$-series of the formal group law. Then there is a fiber sequence for homotopy orbits
    \[(j_\ell^*M)_{hC_\ell} \to M_{hS^1} \xrightarrow{[\ell](u)} M_{hS^1}[2],\]
    and a morphism of fiber sequences
    \begin{equation}\label{eq: diagram fiber sequence gysin}\begin{tikzcd}
    j^*M \arrow[r] \arrow[d]  & {M_{hS^1}} \arrow[d, "\Id"] \arrow[r, "u"]  & M_{hS^1}[2]  \arrow[d, "{[\ell](u)/u}"] \\
    (j_\ell^*M)_{hC_\ell} \arrow[r] & {M_{hS^1}} \arrow[r, "{[\ell](u)}"] & M_{hS^1}[2]
    \end{tikzcd}\end{equation}
    that is compatible with inclusions $C_\ell \rightarrow C_{k\ell}$ for all $k\in \bN$.

    In particular, when $[\ell](u)$ is null homotopic inside $\cE$, there is an isomorphism
    \[(j_\ell^*M)_{hC_\ell} \xrightarrow{\sim} M_{hS^1} \oplus M_{hS^1}[1] .\]

    Similarly, there is a fiber sequence for homotopy fixed points
    \[M^{hS^1}[-2] \xrightarrow{[\ell](u)} M^{hS^1} \to (j_\ell^*M)^{hC_\ell}\]
    that is compatible with inclusions $C_\ell \rightarrow C_{k\ell}$. In particular, when $[\ell](u)$ is null homotopic inside $\cE$, the sequence splits.
\end{Prop}
\begin{proof}
    We follow a similar argument as in \cite[Lemma 7.5]{Kaledin}. Consider the inclusion $i_\ell: \boldsymbol{\Lambda}_\ell \to \boldsymbol{\Lambda}, [n] \mapsto [n \ell]$, where $j_\ell$ factors as $\boldsymbol{\Delta}_{C_\ell} \to \boldsymbol{\Lambda}_\ell \to \boldsymbol{\Lambda}$, and abuse the notations by writing $j_\ell: \boldsymbol{\Delta}_{C_\ell} \to \boldsymbol{\Lambda}_\ell$. 
    Consider the commutative diagram of categories
    \[\begin{tikzcd}
    \boldsymbol{\Delta}_{\bZ/\ell} \arrow[r, "{j_\ell}"] \arrow[d, "\pi"] & \boldsymbol{\Lambda}_\ell \arrow[d, "\pi_\ell"] \ar[r, "i_\ell"] & \boldsymbol{\Lambda} \\
    * \arrow[r, "j"] & \boldsymbol{\Lambda} &.
    \end{tikzcd}\]
    By the base change formula \autoref{lem:base-change},
    we know $(j_\ell^*M)_{hC_\ell} = \pi_{!} j_\ell^*i_\ell^*M = {j}^*\pi_{\ell!}i_\ell^*M$. Thus, for an $\ell$-cyclic object $M$, applying the Gysin sequence for the cyclic object $\pi_!M$ induces the sequence
    \[\pi_!(j_\ell^*i_\ell^* M) \to \pi_!(\pi_{\ell!}i_\ell^*M) \xrightarrow{u} \pi_!(\pi_{\ell!}i_\ell^*M)[2] .\]
    We know that $\pi_!(\pi_{\ell!}i_\ell^*M) = \pi_!(i_\ell^*M)$. Since $\pi: \boldsymbol{\widetilde{\Lambda}}_\ell \to \boldsymbol{\widetilde{\Lambda}}$ is induced by the $\ell$-fold cover of $S^1$, we have $\pi^*(u) = [\ell](u)$ in $H^2(\boldsymbol{\widetilde{\Lambda}}_\ell, E)$, where $[\ell](u)$ represents Chern class after pull-back along the $\ell$-fold cover. Hence, we get the Gysin sequence for $\ell$-cyclic objects
    \[\pi_!(j_\ell^*i_\ell^*M) \to \pi_!(i_\ell^*M) \xrightarrow{[\ell](u)} \pi_!(i_\ell^*M)[2].\]
    
    Finally, we only need to explain why the above sequence is equivalent to the Gysin sequence in the statement. In the end, recall that we have the inclusion $i_\ell: \boldsymbol{\Lambda}_\ell \to \boldsymbol{\Lambda}, [n] \mapsto [n \ell]$, following \cite[Lemma 7.1]{Kaledin}, by the base change formula \autoref{lem:base-change} again, we know $\pi_!(i_\ell^*M ) \simeq \pi_!M$. Thus we get
    \[(j_\ell^*M)_{hC_\ell} \to M_{hS^1} \xrightarrow{[\ell] (u)} M_{hS^1}[2].\]
    Then we can conclude the Gysin sequence mod $\ell$. When $[\ell](u) $ is null homotopic inside $\cE$, it is clear that the mod $\ell$ Gysin sequence splits, which proves the last claim.\end{proof}
    \begin{RMK} 
    We can also achieve the result by the following explicit chain level computation valid over a discrete ring\footnote{We learn this construction from \cite[Proposition 5.3]{Chen1}.}. We denote $S^1_\ell$ the $C_\ell$-invariant cellular decomposition of $S^1$: it consists of $\ell$ 0-cells, and $\ell$ 1-cells cyclic invariant on $S^1$. The corresponding cellular chain $S_\ell=C_{-*}^{cell}(S^1_\ell)$ is represented by $\Bbbk[\tau,\sigma]$ with relations $\tau^\ell=1$, $d\sigma=\tau-1$ and $\sigma^2=0$. The $\ell$-fold covering map $f:S^1_\ell\rightarrow S^1_1, z\mapsto z^\ell$ is cellular, which induces a cochain map $f_*: S_\ell\rightarrow S_1$, where $S_1$ is the algebra of dual numbers $\Bbbk[\varepsilon]$. By the independence of cellular decomposition for computing homologies, we see that $S_\ell$ and $S_1$ are quasi-isomorphic. 
    Moreover, we have $f_*([\ell]_\tau\sigma)=\ell\varepsilon$, where $[\ell]_\tau\coloneqq 1+\tau+\cdots+\tau^{\ell-1}$, which gives the vertical $*\ell$ in \autoref{eq: diagram fiber sequence gysin} after tensoring with $M$ over $\boldsymbol{\Lambda}$. Actually, both $[\ell]_\tau\sigma$ and $\varepsilon$ represent the fundamental class $[S^1]$.  

    For the horizontal $*\ell$ in \autoref{eq: diagram fiber sequence gysin}, let us first note that the Gysin map $u:\Bbbk\rightarrow \Bbbk[2]$ denotes the following chain map using $C_1$ resolution (i.e.~now we have $\tau=1$ and $\sigma=\varepsilon$): 
    \[\begin{tikzcd}
{\Bbbk[\tau,\sigma]} \arrow[r, "{[\ell]_\tau\sigma}"] \arrow[d] & {\Bbbk[\tau,\sigma]} \arrow[r, "{[\ell]_\tau\sigma}"] \arrow[d] & {\Bbbk[\tau,\sigma]} \arrow[r, "{[\ell]_\tau\sigma}"] \arrow[d, "0"] & {\Bbbk[\tau,\sigma]} \arrow[r, "{0}"] \arrow[d] & 0 \arrow[d] \\
{\Bbbk[\tau,\sigma]} \arrow[r, "{[\ell]_\tau\sigma}"]           & {\Bbbk[\tau,\sigma]} \arrow[r, "{0}"]                             & 0 \arrow[r]                                                          & 0 \arrow[r]                                                       & 0  ,        
\end{tikzcd}\]
However if we use $C_\ell$-resolution for $\ell \geq 2$, the above morphism gives us $\ell u$ since $[\ell]_\tau\sigma$ becomes $*\ell$ after tensoring with $M$ (over ${\boldsymbol{\Lambda}}$). To see the fiber $M_{hC_\ell}$, we still use the resolution and base change along $\Bbbk[C_\ell]\simeq \Bbbk[\tau]\to \Bbbk[\tau,\sigma]$ defines an isomorphism between resolutions (over $\Bbbk[\tau]$) of the trivial $\Bbbk[\tau]$-modules $\Bbbk :\Bbbk[\tau,\sigma]^\bullet\rightarrow\Bbbk$ and 
\[\cdots \rightarrow \Bbbk[\tau]\xrightarrow{[\ell]_\tau}  \Bbbk[\tau]\xrightarrow{\tau-1} \Bbbk[\tau]\xrightarrow{\operatorname{aug}} \Bbbk.\]
\end{RMK}

 \subsection{Product formula for cyclic spectra}\label{section: Künneth in general}
Given two cyclic objects $M: \boldsymbol{\Lambda} \rightarrow \cE$ and $N: \boldsymbol{\Lambda} \rightarrow \cE$ where $\cE$ is a symmetric monoidal category, one can define the product as 
\[M \otimes N: \boldsymbol{\Lambda}\xrightarrow{\Delta}\boldsymbol{\Lambda}\times \boldsymbol{\Lambda} \xrightarrow{M \times N} \cE \times \cE \xrightarrow{\otimes_\cE} \cE.\]
Classically, the diagonal map is realized via the cyclic shuffle operation \cite[Chapter 4]{LodayCyclic} \cite{Getzler-Jones}.

We then prove the K\"unneth formula for the homotopy orbits of cyclic objects. See for example \cite[Theorem 4.3.11]{LodayCyclic} or \cite[Theorem 7.29 \& 7.43]{MNN_nilpotence_equivariant_homotopy}. This is the following proposition:

\begin{Prop}\label{prop:Künneth-cyclic}
    Let $M$ and $N$ be cyclic objects over a complex oriented ring spectrum $E$. Then there are fiber sequences (where all the tensor products mean derived tensor products)
    \[(M \otimes N)_{hS^1} \xrightarrow{\Delta} M_{hS^1} \otimes N_{hS^1} \xrightarrow{u \otimes 1 - 1 \otimes u} M_{hS^1} \otimes N_{hS^1}[2].\]
  \[  M^{hS^1} \widehat\otimes N^{hS^1}[-2] \xrightarrow{u \otimes 1 - 1 \otimes u} M^{hS^1} \widehat\otimes N^{hS^1}  \xrightarrow{\Delta}  (M \otimes N)^{hS^1}.\]
    In particular, there are the K\"unneth formulas 
    \[(M \otimes N)_{hS^1} \simeq M_{hS^1} \otimes_{E^{hS^1}} N_{hS^1}[1],\qquad (M \otimes N)^{hS^1} \simeq M^{hS^1} \widehat{\otimes}_{E^{hS^1}} N^{hS^1}.\]
\end{Prop}
\begin{proof}We prove the fiber sequence for homotopy orbits first.

Consider the equivalence $\cE^{BS^1} \simeq \cE \otimes \An^{BS^1}$. We can identify the cartesian cyclic objects $M$ and $N \in \cE^{BS^1}$ as $E \boxtimes M$ and $E \boxtimes N \in \cE \otimes \An^{BS^1}$. Under the straightening/unstraightening equivalence $\An^{BS^1} \simeq \An_{/BS^1}$, $M$ and $N$ correspond to the cartesian fibrations $f_M: M_{hS^1} \to BS^1$, $f_N: N_{hS^1} \to BS^1$ and the product fibration is $f: M_{hS^1} \times N_{hS^1} \to BS^1 \times BS^1$. Here,
    $j^*M = M$ and $j^*N$ are the homotopy fibers.
    Then, by \autoref{lem:thom-pull-back} and \autoref{lem:Künneth}, the Thom spectrum of the pullback bundle $\cO_{\bC P^\infty}(1) \boxtimes \cO_{\bC P^\infty}(-1) \to \bC P^\infty \times \bC P^\infty$ is
    \[\Thom(f^*\cO_{\bC P^\infty}(1) \boxtimes \cO_{\bC P^\infty}(-1)) = \operatorname{cofib}(\Sigma^\infty_+(M \times N)_{hS^1} \to \Sigma^\infty_+M_{hS^1} \otimes \Sigma^\infty_+N_{hS^1}).\]
    Therefore, since $\cE$ is a complex oriented symmetric monoidal category, the Thom isomorphism implies the Gysin sequence
    \[E \otimes \Sigma^\infty_+(M \times N)_{hS^1} \to E \otimes \Sigma^\infty_+M_{hS^1} \otimes \Sigma^\infty_+N_{hS^1} \xrightarrow{u \otimes 1 - 1 \otimes u} E \otimes \Sigma^\infty_+M_{hS^1} \otimes \Sigma^\infty_+N_{hS^1}[2].\]
    The result follows from the equivalence $\cE_{/BS^1} \simeq \cE \otimes \An_{/BS^1}$. The K\"unneth formula follows from the fact that the cofiber of $M_{hS^1} \otimes N_{hS^1} [-2]\rightarrow M_{hS^1} \otimes N_{hS^1}$ is $M_{hS^1} \otimes_{E^{hS^1}} N_{hS^1}$.

    For the fixed points statement, we first obtain a similar fiber sequence:
    \[  M^{hS^1}  \widehat\otimes N^{hS^1}[-2] \xrightarrow{u \otimes 1 - 1 \otimes u} M^{hS^1}  \widehat\otimes N^{hS^1}  \xrightarrow{\Delta}  M^{hS^1} \widehat\otimes_{E^{hS^1}} N^{hS^1}.\]

    We then $u$-adic completion of the cofiber $K\coloneqq \operatorname{cofib}(M^{hS^1} \widehat\otimes N^{hS^1}[-2] \rightarrow M^{hS^1} \widehat\otimes N^{hS^1}$) (which is equivalent to $M^{hS^1} \widehat\otimes N^{hS^1}[-2] \rightarrow M^{hS^1} \widehat\otimes N^{hS^1}$) is $(M \otimes N)^{hS^1}$. 

In fact, by the Gysin sequence \autoref{prop:Gysin}, we have that
\[j^*K\simeq K/(u)\simeq M^{hS^1}/(u) \otimes N^{hS^1}/(u)\simeq M\otimes N \simeq j^*(M\otimes N )\simeq (M \otimes N)^{hS^1}/(u).\]
Consider the homotopy groups. By the derived Nakayama lemma \cite[Tag 0G1U]{StacksProject}, we have the $u$-adic completion of $K$, i.e. $M^{hS^1} \widehat\otimes_{E^{hS^1}} N^{hS^1}$, is equivalent to $(M \otimes N)^{hS^1}$.
\end{proof}

\begin{RMK}    
    One can also show the K\"unneth fiber sequence in the following way. Consider the fibration $M_{hS^1} \times N_{hS^1} \rightarrow BS^1 \times BS^1$. The Atiyah–Hirzebruch spectral sequence implies that $H_*(j^*M) \otimes H_*(E)[u^{-1}] \Rightarrow H_*(M_{hS^1})$ and $H_*(j^*N) \otimes H_*(E)[u^{-1}] \Rightarrow H_*(N_{hS^1})$. Consider the fiber sequence 
    \begin{equation}E[u] \xrightarrow{\Delta} E[u_1^{-1}]] \otimes_E E[u_2^{-1}]] \xrightarrow{u_1^{-1} \otimes 1 - 1 \otimes u_2^{-1}} E[u_1^{-1}]] \otimes_E E[u_2^{-1}]] [2].\end{equation}
    Using functoriality of the spectral sequences, we then obtain the sequence. Also, one can argue directly over a discrete ring $\Bbbk$ as in \cite[Theorem 4.3.11]{LodayCyclic}.
\end{RMK}

\section{Circle Action on Traces}\label{sec:trace}

In this section, we recall the notion of traces of the identity for dualizable objects and the circle action on the trace of the identity, following \cite{Hoyois-Scherotzke-Sibilla,Ayala-Francis-FH-dim1}. This gives an appropriate model for the Hochschild homology with circle action for a dualizable presentable category that are not necessarily compactly generated.

\subsection{Duality and traces}
Recall that an object $X$ in a symmetric monoidal category $(\cC, \otimes, 1)$ is called dualizable if there is another object $X^\vee$ and morphisms
\[\eta_X: 1_\cC \to X^\vee \otimes X, \; \epsilon_X: X \otimes X^\vee \to 1_\cC\]
that satisfy the triangle identities \cite[Section 4.6]{HA} \cite[Section 2.1]{Hoyois-Scherotzke-Sibilla}. One can easily show that when the dual exists, it is unique.
For a morphism $f: X \to Y$ between dualizable objects, we identify its image under $\Hom(X, Y) = \Hom(1_\cC, X^\vee \otimes Y) = \Hom(Y^\vee, X^\vee)$ as the morphism $f^\vee: Y^\vee \to X^\vee$.

For this section, we let $(\cE, \otimes, E)$ be a stable presentable symmetric monoidal category and $(\cC, \otimes, 1)$ is a symmetric monoidal stable category linear over $\cE$.
In fact, we will only consider the case when $\cE$ is a compactly generated symmetric monoidal stable category and $\cC = \Mod_\cE(\PrLst)$ is the symmetric monoidal category of stable presentable categories and colimit preserving functors linear over $\cE$. In this case, $1_\cC = \cE$ and $\End_\cC(1_\cC) = \End_\cC(\cE) = \cE$.

\begin{Def}
For a dualizable object $X$ and an endomorphism $f$, 
we define its trace $\Tr(f) \in \cE$ to be the composition
\[\cE \xrightarrow{\eta_X} X \otimes X^\vee \xrightarrow{f \otimes \,\Id_X} X \otimes X^\vee = X^\vee \otimes X \xrightarrow{\epsilon_X} \cE.\]
We call $\Tr(f)$ the Hochschild homology of $X$ with value in $f$ and denote it by $\HH_*(X, f)$. In particular, we denote $\Tr(\Id_X)$ by $\HH_*(X)$.
\end{Def}

For a morphism $f: X \rightleftharpoons Y$ that admits a right adjoint $f^r$, there is a map $\HH_*(X) \rightarrow \HH_*(Y)$ which is given by the composition \cite[Lemma~4.1]{Brav-Dyckerhoff2}
\begin{align*}
\HH_*(X) &=\Tr(\Id_X) \to \Tr(\Id_X f^rf) \xrightarrow{\sim} \Tr(f \Id_X f^r) \\
&\xrightarrow{\sim} \Tr(\Id_Y ff^r) \to \Tr(\Id_Y) = \HH_*(Y).
\end{align*}
Therefore, Hochschild homology is functorial under right dualizable morphisms.

Consider an idempotent $e: X \to Y \to X$ which can be written as the composition of some inclusion $Y \to X$ and some retraction $X \to Y$ such that both morphisms admit right adjoints. We have the following lemma.

\begin{Lemma}[{\cite[Lemma~2.1]{Hochschild-Kuo-Shende-Zhang}}]\label{lemma: formula for dualizaing data}
Let $(X,X^\vee,\epsilon_X,\eta_X)$ be a duality pair in $\cC$, and
let $e: X \rightarrow X$ be an idempotent which can be written as $X \xrightarrow{r} Y \xrightarrow{i} X$
for some inclusion $i$ and some retraction $r$. 
Assume that the dual idempotent $e^\vee: X^\vee \rightarrow X^\vee$ also splits to $X^\vee \xrightarrow{s} Z \xrightarrow{j} X^\vee$. 
Then the pair
\[\eta_Y \coloneqq (s \otimes r) \circ \eta_X : 1_{\cC} \rightarrow Z \otimes Y, \;\; \epsilon_Y \coloneqq \epsilon_X \circ (i \otimes j) : Y \otimes Z \rightarrow 1_{\cC}\]
exhibits $Z$ as the dual of $Y$. 
\end{Lemma}

We note that Hochschild homology is a localizing invariant. Namely, we have the following formula:

\begin{Prop}[{\cite[Theorem 3.4]{Hoyois-Scherotzke-Sibilla}}]\label{prop:hochschild-fiber-seq}
Let $X, Y$ and $Z$ be dualizable objects in $\cC$ such that there is a localization sequence $X \to Y \to Z$ in $\cC$.
Then there is a cofiber sequence
\[\HH_*(X) \to \HH_*(Y) \to \HH_*(Z).\]
\end{Prop}

Finally, we note that for any object $X$ in $\cC$, we can also define the Hochschild cohomology as follows:

\begin{Def}
For any object $X$ and an endomorphism $f:X\to X$, we define the Hochschild cohomology by $\HH^*(X, f) = \Hom(\Id_X, f)$. In particular, we denote $\HH^*(X, \Id_X) = \Hom(\Id_X, \Id_X)$ by $\HH^*(X)$.
\end{Def}

\subsection{\texorpdfstring{$S^1$-action on traces}{}}
The trace of a dualizable category admits a natural $S^1$-action. We recall the construction of $S^1$-action using factorization homology by Ayala--Francis \cite{Ayala-Francis-FH-dim1}, following earlier works of \cite{Ayala-Francis-FH,AMGR-FH-enriched}\footnote{We would like to thank Nick Rozenblyum for pointing out the reference of the construction, and we would like to thank Nick Rozenblyum and John Francis for explaining to us these works.}.

    We recall the topological models of the para-cyclic categories from factorization homology considered in Ayala--Mazel-Gee--Rozenblyum \cite{AMGR-FH-enriched} and Ayala--Francis \cite{Ayala-Francis-FH-dim1} in the following remarks.

\begin{RMK}\label{rem:cyclic-combin-topol}
    Let $\cD(S^1)$ be the category of framed 1-disk refinements of $S^1$, whose objects are stratifications of $S^1$ by disks and morphisms are maps between stratified spaces \cite[Definition 1.18]{AMGR-FH-enriched}. By \cite[Observation 1.20]{AMGR-FH-enriched} or \cite[Observation 1.3 \& Corollary B.9]{Ayala-Francis-FH-dim1}, there exists an $S^1$-equivariant equivalence between $\infty$-categories (where $sk_0$ denotes the 0-skeleton of the stratification)
    \[\cD(S^1) \simeq \boldsymbol{\Delta}_\circlearrowleft, \; (D \to S^1) \mapsto \pi_0(\exp^{-1}(S^1 \setminus sk_0(D))).\]
\end{RMK}
\begin{RMK}\label{rem:disk-stratum-embed}
    Consider the categories $\cD(S^1)^{\triangleleft}$ and $\cD(S^1)^{\triangleright}$ defined by adjoining an initial or a final object to $\cD(S^1)$. Let $\operatorname{Disk}_{1/S^1}^\text{fr}$ be the category of framed embedded disks into $S^1$, whose objects are disjoint unions of framed embedded disks in $S^1$ and morphisms are framed embeddings \cite[Section 2.2]{Ayala-Francis-FH}. By \cite[Corollary B.9]{Ayala-Francis-FH-dim1}, there is an $S^1$-equivariant equivalence
    \[\cD(S^1)^\triangleleft \simeq \operatorname{Disk}_{1/S^1}^\text{fr}, \;\; \{*\} \mapsto (\varnothing \hookrightarrow S^1),\;\; (D \to S^1) \mapsto (D \setminus sk_0(D) \hookrightarrow S^1).\]
\end{RMK}

\begin{Thm}[{\cite[Theorem 2.6]{Ayala-Francis-FH-dim1}}]\label{thm:AF}
    Let $\cC$ be a symmetric monoidal category enriched over $\cE$ that admits geometric realizations. Let $X$ be a dualizable object in $\cC$. Then there is a canonical functorial $S^1$-action on $\HH_*(X) = \Tr(\Id_X)$ defined by a cyclic object $\Tr(\Id_X, -) : \boldsymbol{\Lambda} \to \cE$ such that
    \begin{align*}
        \Tr(\Id_X, n): \cE \xrightarrow{\eta_X^{\otimes n}} (X^\vee \otimes X)^{\otimes n} \xrightarrow{\tau_{2n}} (X \otimes X^\vee)^{\otimes n} \xrightarrow{\epsilon_X^{\otimes n}} \cE,
    \end{align*}
    where $\tau_{2n}: (X^\vee \otimes X)^{\otimes n} \xrightarrow{\sim} (X \otimes X^\vee)^{\otimes n}$ is the cyclic permutation of order $2n$, such that the cyclic rotation $t_n$ is induced by $\tau_{2n}^2$, and the face maps $d_i$ and degeneracy maps $s_i$ are induced by the triangle identities of the unit and counit.
\end{Thm}
\begin{RMK}
    \label{rem:cyclic-object-map}
    The face morphism $d_i: \Tr(\Id_X, n) \to \Tr(\Id_X, n+1)$ is induced by the identity 
    \[(\epsilon_X \otimes \epsilon_X) \circ (\Id_X \otimes \,\eta_X \circ \Id_{X^\vee}) \xrightarrow{\sim} \epsilon_X.\]
    The degeneration morphism $s_i: \Tr(\Id_X, n+1) \to \Tr(\Id_X, n)$ is induced by the identity 
    \[\eta_X \xrightarrow{\sim} (\Id_X \otimes \,\epsilon_X \circ \Id_{X^\vee}) \circ (\eta_X \circ \,\eta_X).\]
    Finally, the cyclic rotation is induced by $\tau_{2n}^2: (X^\vee \otimes X)^{\otimes n} \to (X^\vee \otimes X)^{\otimes n}$ composed with the units and precomposed with the counits.
\end{RMK}

\begin{proof}[Proof of \autoref{thm:AF}]
    Consider $\cD(S^1)$ and $\operatorname{Disk}_{1/S^1}^\text{fr}$. 
    Following \autoref{rem:disk-stratum-embed}, the factorization homology of $X$ is defined by \cite[Observation 2.3]{Ayala-Francis-FH-dim1}
    \[\HH_*(X) \coloneqq \colim \left( \cD(S^1)^{\triangleleft} \xrightarrow{\sim} \operatorname{Disk}_{1/S^1}^\text{fr} \to \operatorname{Disk}_{1}^\text{fr} \xrightarrow{\inEnd(X)} \cC \right).\]
    Ayala--Francis \cite[Theorem 2.6]{Ayala-Francis-FH-dim1} showed that there exists an $S^1$-invariant trace morphism from the factorization homology to $\cV$, or equivalently, that the colimit of the factorization homology over $\cD(S^1)^\triangleleft$ factors through $\cD(S^1)^{\triangleleft \triangleright}$:
    \[\HH_*(X) \coloneqq \colim \left( \cD(S^1)^{\triangleleft} \xrightarrow{\inEnd(X)} \cC \right) = \colim \left( \cD(S^1)^{\triangleleft \triangleright} \xrightarrow{\inEnd(X)} \cC \right).\]
    In particular, there exists a functor $\Tr(f, -): \cD(S^1)^{\triangleleft \triangleright} \to \cV$ by evaluating in the  category of framed 1-disk refinements after adjoining an initial and a final object $\cD(S^1)^{\triangleleft \triangleright}$. More precisely, since $X$ is dualizable, we have that $\inEnd(X) = X^\vee \otimes X$, it follows from \cite[Equation (10) \& (11)]{Ayala-Francis-FH-dim1} that the evaluation at the map from the initial object to the disk stratification $D \to S^1$ and then to the final object is given by
    \begin{align*}
        \cE \xrightarrow{\eta_X^{\otimes sk_0(D)}} (X^\vee \otimes X)^{\otimes \pi_0(D)} \xrightarrow{\tau_{2n}} (X \otimes X^\vee)^{\otimes \pi_0(S^1 \setminus sk_0(D))} \xrightarrow{\epsilon_X^{\otimes \pi_0(S^1 \setminus sk_0(D))}} \cE,
    \end{align*}
    where $\tau_{2n}: sk_0(D) \xrightarrow{\sim} \pi_0(S^1 \setminus sk_0(D))$ is induced by Poincar\'e duality exchanging the 0-dimensional and 1-dimensional strata in $S^1$. Equivalently, we have a desuspension functor $\Tr(\Id_X, -): \cD(S^1) \to \cV$. For a pair of disk refinements $D \to D' \to S^1$, we have natural morphisms $\Tr(\Id_X, D) \to \Tr(\Id_X, D')$ that are induced by the triangle identities of the unit and counits.
    
    Let $\boldsymbol{\Delta}_\circlearrowleft$ be the paracyclic category. By \autoref{rem:cyclic-combin-topol}, there exists an $S^1$-equivariant equivalence $\cD(S^1) \simeq \boldsymbol{\Delta}_\circlearrowleft$.
    Thus we get an $S^1$-invariant paracyclic object $\Tr(\Id_X, -) : \boldsymbol{\Delta}_\circlearrowleft \to \cV$ such that
    \begin{align*}
        \Tr(\Id_X, n): \cE \xrightarrow{\eta_X^{\otimes n}} (X^\vee \otimes X)^{\otimes n} \xrightarrow{\tau_{2n}} (X \otimes X^\vee)^{\otimes n} \xrightarrow{\epsilon_X^{\otimes n}} \cE,
    \end{align*}
    Since $(\boldsymbol{\Delta}_\circlearrowleft)_{hS^1} = \boldsymbol{\Lambda}$ by  \autoref{prop:paracyclic}, we obtain a cyclic object $\Tr(\Id_X, -) : \boldsymbol{\Lambda} \to \cE$. Thus, we explained how the theorem follows from \cite[Theorem 2.6]{Ayala-Francis-FH-dim1}.
\end{proof}

\begin{Def}
For a dualizable object $X$ in the symmetric monoidal category $\cC$, we define the cyclic homology to be the $S^1$-homotopy orbits of the Hochschild homology $\HH_*(X)_{hS^1}$, and the negative cyclic homology to be the $S^1$-homotopy fixed points of the Hochschild homology $\HH_*(X)^{hS^1}$.
\end{Def}

\begin{RMK}
More generally, by Hoyois--Scherotzke--Sibilla \cite[Proposition 2.16]{Hoyois-Scherotzke-Sibilla}, which uses the cobordism hypothesis \cite{Lurie-cobordism}, the trace functor admits a canonical $S^1$-invariant lift $\Tr: \Aut(\cC) \to \End_\cC(1_\cC)$ with respect to the canonical $S^1$-action on $\Aut(\cC)$. For any single automorphism $f: X \to X$ that is an $S^1$-fixed point in $\Aut(\cC)$, there is an $S^1$-action on the trace $\Tr(f)$ \cite[Proposition 4.3]{Brav-Dyckerhoff2}.
\end{RMK}

\subsection{\texorpdfstring{$C_\ell$-action on traces}{}}
For an automorphism of dualizable objects, we have shown that there is an $S^1$-action on the trace. Let $\ell$ be a natural number. In this section, we study the $C_\ell$-action on the trace by restriction along the standard inclusion $C_\ell \to S^1$.

We can consider the $C_\ell$-actions described by $C_\ell$-category $\boldsymbol{\Delta}_{C_\ell}$. By \autoref{lem:finite-action-1}, we have an immediate corollary:

\begin{Prop}\label{prop:inclusion-in-circle}
    Let $\cC$ be a symmetric monoidal category enriched over $\cE$ that admits geometric realizations. Let $X$ be a dualizable object in $\cC$. Then there is a canonical functorial $C_\ell$-action on $\HH_*(X) = \Tr(\Id_X)$ defined by a $\ell$-cyclic object $\Tr_\ell(\Id_X) : \boldsymbol{\Delta}_{\bZ/\ell} \to \cE$ with
    \begin{align*}
        \Tr_\ell(\Id_X, \ell): \cE \xrightarrow{\eta_X^{\otimes \ell}} (X^\vee \otimes X)^{\otimes \ell} \xrightarrow{\tau_{2\ell}} (X \otimes X^\vee)^{\otimes \ell} \xrightarrow{\epsilon_X^{\otimes \ell}} \cE,
    \end{align*}
    where $\tau_{2\ell}: (X^\vee \otimes X)^{\otimes \ell} \xrightarrow{\sim} (X \otimes X^\vee)^{\otimes \ell}$ is the cyclic permutation of order $2\ell$, such that $t_\ell$ is induced by the cyclic rotation on $(X^\vee \otimes X)^{\otimes \ell}$ of order $\ell$.
\end{Prop}

Given any $S^1$-action, we can naturally obtain another $S^1$-action by considering the $\ell$-fold cover, described by the $\ell$-cyclic categories $\boldsymbol{\Lambda}_\ell$. By \autoref{prop:ellcyclic}, we have the following immediate corollary:

\begin{Coro}
    Let $\cC$ be a symmetric monoidal category enriched over $\cE$ that admits geometric realizations. Let $X$ be a dualizable object in $\cC$. Then for any $\ell \in \bN^\times$, there is a canonical functorial $S^1$-action on $\HH_*(X) = \Tr(\Id_X)$ defined by a $\ell$-cyclic object $\Tr(\Id_X, -) : \boldsymbol{\Lambda}_\ell \to \cE$ such that
    \begin{align*}
        \Tr(\Id_X, n): \cE \xrightarrow{\eta_X^{\otimes n\ell}} (X^\vee \otimes X)^{\otimes n\ell} \xrightarrow{\tau_{2n\ell}} (X \otimes X^\vee)^{\otimes n\ell} \xrightarrow{\epsilon_X^{\otimes n\ell}} \cE.
    \end{align*}    
\end{Coro}

\begin{RMK}
    When $\cE = \Sp$ and $\cC = \PrLst$, we will see that the $\ell$-cyclic traces assemble into a cyclotomic object via the (generalized) Tate diagonal construction \cite[Chapter III.1]{Nikolaus-Scholze-Ontopologicalcyclichomology} \cite[Section 3.1]{AMGR-naive} \cite[Section 4.1--4.3]{AMGR-cyclotomic}. For a prime $p$, that is to say, there is a lax commutative diagram of cyclic objects:
    \[\varphi_p: \Tr(\Id_X) \to \Tr(\Id_X)^{tC_p}:= \mathrm{cofib}(\Tr(\Id_X)_{hC_p} \to \Tr(\Id_X)^{hC_p}).\]
    We will discuss the cyclotomic structure in the next section using the cyclic bar construction.
\end{RMK}

Finally, following the result of Chen \cite{Chen2}, we can define the $p$-th power operation on Hochschild cohomology and the $p$-fold cap product on the $C_p$-invariant Hochschild homology on the level of spectra.

\begin{Prop}[{\cite[Section 2]{Chen2} \cite[Section 3]{Rezchikov}}]\label{prop:power construction}Let $p$ be a prime number.
    (1)~Let $\cC$ be a symmetric monoidal category enriched over \textbf{$\bZ/p\bZ$} that admits geometric realizations. Let $X$ be a dualizable object in $\cC$. Then there is a $p$-fold cap operation 
    \[\cap_p:\HH^*(X) \to \Hom(\HH_*(X)^{tC_p}, \HH_*(X)^{tC_p})\]
    or equivalently a morphism
    \[\Hom(\Id_X, \Id_X) \to (\Hom(\Id_X, \Id_X)^{\otimes p})^{tC_p} \to \Hom(\Tr(\Id_X)^{tC_p}, \Tr(\Id_X)^{tC_p}).\]
    (2)~Let $\cC$ be a symmetric monoidal category enriched over \textbf{$\Sp$} that admits geometric realizations. Let $X$ be a dualizable object in $\cC$. Then there is a $p$-fold cap operation 
    \[\cap_p:\HH^*(X) \to \Hom(\HH_*(X)^{tC_p}, \HH_*(X)^{tC_p})\]
    or equivalently a morphism
    \[\Hom(\Id_X, \Id_X) \to (\Hom(\Id_X, \Id_X)^{\otimes p})^{tC_p} \to \Hom(\Tr(\Id_X)^{tC_p}, \Tr(\Id_X)^{tC_p}).\]
\end{Prop}
\begin{proof}
    (1)~We have a map $T_p: \Hom(\Id_X, \Id_X) \to (\Hom(\Id_X, \Id_X)^{\otimes p})^{hC_p}$ defined by the $p$-th power operation, which is only additive after inverting $u \in H^2(\widetilde{\boldsymbol{\Delta}}_{\bZ/p}, \bZ/p\bZ)$ \cite[Lemma 2.5]{SeidelWilkins}. By \autoref{prop:Gysin-finite-group}, it induces a Frobenius-linear map $\Hom(\Id_X, \Id_X) \to (\Hom(\Id_X, \Id_X)^{\otimes p})^{tC_p}$. Consider the cap product on the $\boldsymbol{\widetilde{\Lambda}}_p$-cyclic objects
    \[\Hom(\Id_X, \Id_X)^{\otimes p} \to \Hom(\Tr_p(\Id_X), \Tr_p(\Id_X)).\]
    By taking the Tate fixed point and the K\"unneth morphism (see for example \cite{BlumbergMandell}), we obtain
    \begin{align*}(\Hom(\Id_X, \Id_X)^{\otimes p})^{tC_p} &\to \Hom(\Tr_p(\Id_X), \Tr_p(\Id_X))^{tC_p} \\
    &\to \Hom(\Tr(\Id_X)^{tC_p}, \Tr(\Id_X)^{tC_p}).
    \end{align*}
    
    (2)~First, we define the map $T_p: \Hom(\Id_X, \Id_X) \to (\Hom(\Id_X, \Id_X)^{\otimes p})^{tC_p}$ from the Tate diagonal morphism on spectra \cite[Section III.1]{Nikolaus-Scholze-Ontopologicalcyclichomology}. Then, consider the cap product on the $\boldsymbol{\widetilde{\Lambda}}_p$-cyclic objects
    \[\Hom(\Id_X, \Id_X)^{\otimes p} \to \Hom(\Tr_p(\Id_X), \Tr_p(\Id_X)).\]
    By taking the Tate fixed point and the K\"unneth morphism \cite{BlumbergMandell}, we obtain
    \begin{align*}(\Hom(\Id_X, \Id_X)^{\otimes p})^{tC_p} &\to \Hom(\Tr_p(\Id_X), \Tr_p(\Id_X))^{tC_p} \\
    &\to \Hom(\Tr(\Id_X)^{tC_p}, \Tr(\Id_X)^{tC_p}).
    \end{align*}
    This completes the proof.
\end{proof}
\begin{RMK}
    We remark that it is shown by \cite[Theorem 1.15]{Rezchikov} plus \cite[Section III.6]{Nikolaus-Scholze-Ontopologicalcyclichomology} that the above constructions over $\mathbb{S}$ and $\bZ/p\bZ$ are related by base change. In other words, the constructions in statement (1) and (2) are related by base change.
\end{RMK}
\begin{RMK}
    It is explained in \cite[Section 4]{Rezchikov} that we have a commutative diagram of the ordinary cap product and $p$-fold cap product
    \[\begin{tikzcd}
    \HH^*(X) \otimes \HH_*(X) \ar[r,"\cap=\cap_1"] \ar[d,""] & \HH_*(X) \ar[d] \\
    (\HH^*(X)^{\otimes p})^{tC_p} \otimes \HH_*(X)^{tC_p} \ar[r,"\cap_p"] & \HH_*(X)^{tC_p}.
    \end{tikzcd}\]
    In fact, the formula can be extracted from \cite[Theorem 1.5]{Rezchikov} and the maps $\Phi^{\bZ/p}(N^{\bZ/p}\HH^*(X)) \to N^{\bZ/p}\HH^*(X)^{tC_p}$ and $\Phi^{\bZ/p}HH_*(X) \to \HH_*(X)^{tC_p}$ \cite[Section III.6]{Nikolaus-Scholze-Ontopologicalcyclichomology}.
\end{RMK}
\begin{RMK}
    It is explained in \cite[Section IV.1]{Nikolaus-Scholze-Ontopologicalcyclichomology} that the $p$-th power construction is closely related to the Steenrod operation
    \begin{align*}
    H^*(M; \bZ/p\bZ) \to H^*(M; \bZ/p\bZ) \otimes H_{C_p}^*(pt; \bZ/p\bZ), \quad \beta \mapsto \mathrm{St}(\beta).
    \end{align*}
    We will see later in \autoref{sec:cat-of-sheaf} how to reconstruct the Steenrod operation from Hochschild cohomology of the category of sheaves. In the case of Fukaya categories, it is also known that the construction  gives the quantum Steenrod operation \cite{Chen2}.
\end{RMK}

\subsection{Gysin sequence and Künneth formula}
In this section, we adapt the result of \autoref{section: Künneth in general} and \autoref{section: Gysin in general} to the circle action on trace.

First, we consider the K\"unneth formula for the Hochschild homology with $S^1$-actions. See for example \cite[Lemma 6.3]{Brav-Dyckerhoff2}.

\begin{Prop}\label{cor:Künneth}
    Let $\cC$ be a symmetric monoidal category enriched over $\cE$ that admits geometric realizations. Let $X$ and $Y$ be dualizable objects in $\cC$. Then there is a canonical isomorphism $\HH_*(X \otimes Y) = \HH_*(X) \otimes \HH_*(Y)$ given by isomorphism of cyclic objects
    \[\Tr(\Id_{X \otimes Y}, -) = \Tr(\Id_X, -) \otimes \Tr(\Id_Y, -).\]
    In particular, when $E$ is a complex oriented ring spectrum, there is a fiber sequence
    \[\HC_*(X \otimes Y) \to \HC_*(X) \otimes \HC_*(Y) \xrightarrow{u \otimes 1 - 1 \otimes u} \HC_*(X) \otimes \HC_*(Y)[2],\]
    and there is the K\"unneth formula $\HC_*(X \otimes Y) \simeq \HC_*(X) \otimes_{E^{hS^1}} \HC_*(Y)[1]$.
\end{Prop}
\begin{proof}
By \autoref{prop:Künneth-cyclic}, this follows immediately from the fact that for dualizable objects $X$ and $Y$, the tensor product $X \otimes Y$ is also dualizable, the unit is $\eta_{X \otimes Y} = \eta_X \otimes \eta_Y$ and the counit is $\epsilon_{X \otimes Y} = \epsilon_X \otimes \epsilon_Y$. 
\end{proof}

\begin{Prop}
    Let $X$ be a dualizable object over a complex oriented symmetric monoidal stable category $\cE$. There is a $C_\ell$-Gysin sequence
    \[\Tr_\ell^E(\Id_X)_{hC_\ell} \to \Tr^E(\Id_X)_{hS^1} \xrightarrow{[\ell] (u)} \Tr^E(\Id_X)_{hS^1}[2].\]
    When $[\ell](u)$ is null homotopic in $\cE$, we have an isomorphism
    \[ \Tr_\ell^E(\Id_X)_{hC_\ell} \simeq \Tr^E(\Id_X)_{hS^1} \oplus \Tr^E(\Id_X)_{hS^1}[1].\]
    In particular, when $E = \bZ$ and $\Tr^\bZ(\Id_X)_{hS^1}$ has no $\ell$-torsion, then for a field $\Bbbk$ of characteristic $char\,\Bbbk\,|\,\ell$, there is an isomorphism
    \[ \Tr_\ell^\Bbbk(\Id_X)_{hC_\ell} \simeq (\Tr^\bZ(\Id_X)_{hS^1} \oplus \Tr^\bZ(\Id_X)_{hS^1}[1] ) \otimes_\bZ \Bbbk.\]
\end{Prop}

\begin{proof}We apply \autoref{prop:Gysin-finite-group} to $\Tr$.
\end{proof}

\subsection{Relation with cyclic bar constructions}
There are several definitions of coherent $S^1$-actions on Hochschild homology in the literature. Using the cobordism hypothesis, all of them should be identical. However, some of the identifications can also be shown directly. 

Using the enriched factorization homology, Ayala--Mazel-Gee--Rozenblyum gave another combinatorial construction of the circle action on Hochschild homology. The following theorem can be extracted from \cite{Hoyois}, \cite[Definition III.2.3]{Nikolaus-Scholze-Ontopologicalcyclichomology}, or \cite[Definition 3.12, Remark 1.21 \& Remark 2.26]{AMGR-FH-enriched}.

\begin{Thm}[\cite{Nikolaus-Scholze-Ontopologicalcyclichomology,AMGR-FH-enriched}]\label{thm:AMGR}
    Let $X$ be a compactly generated stable $\infty$-category linear over $\cE$\footnote{In fact, Ayala--Mazel-Gee--Rozenblyum \cite{AMGR-FH-enriched} considered $\infty$-categories enriched over $\cE$. However, when $X$ is stable, being enriched over $\cE$ is equivalent to being linear over $\cE$. This is for instance the setting in \cite{Hoyois-Scherotzke-Sibilla}.}. Then there is a functorial $S^1$-action on $\HH_*(X) \coloneqq \Tr(\Id_X)$ defined by a cyclic object $\HH_*(\Id_X, -): \boldsymbol{\Lambda} \to \cE$ with
    \begin{align*}
        \HH_*(\Id_X, n) = \mathrm{Bar}_n(X, \Id_X) =\bigoplus_{x_0, \dots, x_{n-1} \in X^c} ~\bigotimes_{1\leq i\leq n} \HOM(x_{i-1}, x_i).
    \end{align*}
\end{Thm}

Finally, the construction of Ayala--Francis using traces \cite{Ayala-Francis-FH-dim1} is compatible with Ayala--Mazel-Gee--Rozenblyum using bar constructions \cite{AMGR-FH-enriched}, as shown in \cite[Proposition 4.24]{Hoyois-Scherotzke-Sibilla}, \cite[Proposition B.12]{Ayala-Francis-FH-dim1} and \cite[Remark 2.3 \& Remark 4.14]{AMGR-FH-enriched}, whose proof we now recall. 

\begin{Prop}[{\cite[Proposition 4.24]{Hoyois-Scherotzke-Sibilla}, \cite[Proposition B.12]{Ayala-Francis-FH-dim1}}]\label{prop:compare-fh}
    Let $X$ be a compactly generated stable $\infty$-category over $\cE$. Then the functorial $S^1$-actions on $\HH_*(X)$ in \autoref{thm:AF} and \autoref{thm:AMGR} agree. In other words, there is an $S^1$-equivariant isomorphism
    \[\Tr_\bullet(\Id_X) \simeq |\mathrm{Bar}_\bullet(X, \Id_X)| =\operatornamewithlimits{colim}_{n\in \Delta} \bigoplus_{x_0, \dots, x_{n-1}} \bigotimes_{1\leq i\leq n} \HOM(x_{i-1}, x_i).\]
\end{Prop}
\begin{proof}
    First, following \cite[Proposition 4.24]{Hoyois-Scherotzke-Sibilla}, we can show that for a compactly generated $\infty$-category $X$ over $\cE$, the evaluation is given by the bar construction, so that for any bimodule $M$, which by definition is an object in $\Fun_\cE(X^{c,op} \otimes X^{c}, \cE)$,
    \[\epsilon_X(M) = \mathrm{Bar}(X, M) = \colim_{n\in \Delta}\bigoplus_{x_0,\dots,x_{n-1}}\Big(\bigotimes_{1\leq i\leq n-1}\Hom(x_{i-1}, x_i)\Big) \otimes M(x_{n-1}, x_0).\]
    In fact, the identity $\epsilon_X(-) = \mathrm{Bar}(X, -)$ follows from $\Id_X = (\Id_X \otimes \mathrm{Bar}(X, -)) \circ (\eta_X \otimes \Id_X)$. The latter is because the natural map
    \[\mathrm{Bar}(X, \Hom(x, -) \otimes N(-)) \to N(x)\]
    defines an augmented simplicial object with an extra degeneracy map sending $N(x) \to \Hom(x, x) \otimes N(x)$, $\Hom(x, a) \otimes N(a) \to \Hom(x, a) \otimes \Hom(x, x) \otimes N(a)$ and thus is contractible.
    
    Let $M$ be the diagonal bimodule $M(x, y) = \Hom(x, y)$ (in other words, it corresponds to the unit $\eta_X$). We can rewrite the trace using the cyclic bar construction:
    \[\Tr(\Id_X, r) \simeq |\mathrm{Bar}_\bullet(X, \Id_X, \dots, \Id_X)| = \colim_{n\in \Delta}\bigoplus_{x_0,\dots,x_{rn-1}}\bigotimes_{1\leq i\leq rn}\Hom(x_{i-1}, x_i).\]
    where the face maps are given by $n$-fold compositions. Therefore, we claim that there is a natural commutative diagram between the functors $\boldsymbol{\Delta}_{\circlearrowleft} \to \cE$ 
    \begin{align*}
    \HH_*(\Id_X, r)
    &= \bigoplus_{x_0, x_1, \dots, x_{r-1} \in X^c} \bigotimes_{1\leq i\leq r} \HOM(x_{i-1}, x_i)\\
    &\to \colim_{n\in \Delta}\bigoplus_{x_0,\dots,x_{rn-1}\in X^c}\bigotimes_{1\leq i\leq rn}\Hom(x_{i-1}, x_i)= \Tr(\Id_X, r) .
    \end{align*}
    While the two cyclic objects are not isomorphic, we will show that the corresponding cartesian cyclic objects are, and for that, it suffices to check the isomorphisms at the point as in \cite[Proposition 4.24]{Hoyois-Scherotzke-Sibilla}, or equivalently, the pull-back along $\boldsymbol\Delta \to \boldsymbol{\Lambda}$ \cite[Proposition 1.1]{Hoyois}:
    \[\colim_{n\in \Delta}\bigoplus_{x_0,\dots,x_{n-1}}\bigoplus_{1\leq i\leq n}\Hom(x_{i-1},x_i) \xrightarrow{\sim} \colim_{r\in\bN}\colim_{n\in\Delta}\bigoplus_{x_0,\dots,x_{rn-1}}\bigotimes_{1\leq i\leq rn}\Hom(x_{i-1}, x_i).\]
    
    Finally, we explain the claim we made on the commutative diagram of paracyclic objects via factorization homology. As explained in the proof of \cite[Proposition B.12]{Ayala-Francis-FH-dim1}, by writing
    \[\Tr(\Id_X, r) \simeq |\mathrm{Bar}_\bullet(X, \Id_X, \dots, \Id_X)| = \colim_{n\in \Delta}\bigoplus_{x_0,\dots,x_{rn-1} }\bigotimes_{1\leq i\leq rn}\Hom(x_{i-1}, x_i),\]
    one deduces that there is a commutative diagram of functors from the (disk model of) paracyclic category $\cD(S^1)$ into $\cE$, where the top horizontal functor computes $\HH_*(\Id_X)$ in  \autoref{thm:AMGR} and the bottom horizontal functor computes $\Tr(\Id_X)$ in \autoref{thm:AF}:
    \[\begin{tikzcd}
    \cD(S^1) \ar[r] \ar[equal,d] & \cD_{/S^1} \ar[r] & \cD \ar[r, "\mathrm{Bar}_\bullet(X)"] & \cE \ar[equal,d] \\
    \cD(S^1) \ar[r] & \mathrm{Disk}^\mathrm{fr}_{1/S^1} \ar[r] & \mathrm{Disk}^\mathrm{fr}_1 \ar[r, "\End(X)"] & \cE.
    \end{tikzcd}\]
    For that claim, it suffices to show that the natural functor $\cD(S^1) \to \mathrm{Disk}^\mathrm{fr}_{1/S^1}$ that sends stratifications of $S^1$ by disks to framed embedded disks in $S^1$ is final. This follows from the proof of \cite[Proposition B.12]{Ayala-Francis-FH-dim1} using \cite[Lemma 2.3.2]{Ayala-Francis-PD} and \cite[Observation B.8]{Ayala-Francis-FH}.
\end{proof}

\begin{RMK}
     \autoref{prop:compare-fh} shows that the cyclic trace of dualizable categories defined in  \autoref{thm:AF} \cite{Ayala-Francis-FH-dim1} is the unique extension of the Hochschild homology of compactly generated categories in \autoref{thm:AMGR} \cite{AMGR-FH-enriched,Nikolaus-Scholze-Ontopologicalcyclichomology} to dualizable categories in the sense of Efimov \cite[Theorem 0.1]{Efimov-K-theory}.
\end{RMK}

Finally, the above construction \cite{AMGR-FH-enriched} and Nikolaus--Scholze \cite{Nikolaus-Scholze-Ontopologicalcyclichomology} can be identified by writing down both constructions explicitly (see \cite[Section 0.3]{AMGR-FH-enriched}). The comparison between this and the classical mixed complex construction by Connes and Tsygan is shown by Hoyois \cite{Hoyois}. A comparison between the trace and the classical bar construction of Hochschild homology, on the other hand, has also appeared in \cite[Proposition 4.24]{Hoyois-Scherotzke-Sibilla}, but they did not compare the $S^1$-actions. In fact, the $S^1$-action defined in \cite[Proposition 2.16]{Hoyois-Scherotzke-Sibilla} invokes the cobordism hypothesis \cite{Lurie-cobordism}.

There has also been a $C_\ell$-action on Hochschild homology in the literature using the mod $\ell$ bar construction by Chen \cite{Chen1}. Namely, Chen constructed a finite cyclic category $\boldsymbol{\Lambda}_{\bZ/\ell}$ (which is denoted by ${}_\ell\widetilde{\boldsymbol{\Lambda}}$ in \cite{Chen1}) whose objects are cyclic sets with $\ell$ marked points and morphisms are maps that send marked points to marked points. For a cyclic set with $\ell$ marked points and $n_i$ points in between the $i$ and $(i+1)$-th marked point, we denote it by $[n_1, \dots, n_\ell]$. There is an obvious factorization
\[\boldsymbol{\Delta}_{\bZ/\ell} \to \boldsymbol{\Lambda}_{\bZ/\ell} \to \boldsymbol{\Lambda}\]

\begin{Prop}[{Chen \cite[Proposition 3.19]{Chen1}}]\label{prop:Chen}
    Let $X$ be a compactly generated $\infty$-category over $\cE$. Then there is a functorial $C_\ell$-action on $\HH_*(X) \coloneqq \Tr(\Id_X)$ defined by a cyclic object $\HH_*(\Id_X, -): \boldsymbol{\Lambda}_{\bZ/\ell} \to \cE$ with
    \begin{align*}
        \HH_*(\Id_X, [n_1, \dots, n_\ell]) = \bigoplus_{x_0, x_1, \dots, x_{n_1+\dots+n_\ell+\ell-1} \in X^c} \bigotimes_{1\leq i\leq n} \HOM(x_{i-1}, x_i).
    \end{align*}
\end{Prop}

Similar to \autoref{lem:finite-action-1} in our setting, Chen also showed that the inclusion $j: \boldsymbol{\Lambda}_{\bZ/\ell} \to \boldsymbol{\Lambda}$ is induced by the standard inclusion $C_\ell \hookrightarrow S^1$.

\begin{Lemma}[{\cite[Lemma 5.1]{Chen1}}]\label{lem:finite-action-2}
Consider the natural inclusion $j: \boldsymbol{\Lambda}_{\bZ/\ell} \to \boldsymbol{\Lambda}$
and the standard inclusion $j: C_\ell \to S^1$. There is a commutative diagram of categories
\[\xymatrix{
\boldsymbol{\Lambda}_{\bZ/\ell} \ar[r] \ar[d] & \boldsymbol{\Lambda} \ar[d] \\
BC_\ell \ar[r] & BS^1.
}\]
\end{Lemma}

Then it follows from \autoref{lem:finite-action-1}, ~\ref{lem:finite-action-2} and \autoref{prop:compare-fh} that the finite cyclic homology defined using traces is isomorphic to the finite cyclic homology defined using bar construction.

\begin{Prop}\label{prop:compare-fh-finite}
    Let $X$ be a compactly generated $\infty$-category over $\cE$. Then the functorial $C_\ell$-actions on $\HH_*(X)$ in \autoref{prop:inclusion-in-circle} and \autoref{prop:Chen} agree.
\end{Prop}
\begin{proof}
    The $C_\ell$-actions on $\Tr(\Id_X)$ in \autoref{prop:inclusion-in-circle} are induced by the pull-back of the $S^1$-action in \autoref{thm:AF} along the map $\boldsymbol{\Delta}_{\bZ/\ell} \to \Lambda$, and the $C_\ell$-actions on $\HH_*(\Id_X)$ in \autoref{prop:Chen} are induced by the pull-back of the $S^1$-action in \autoref{thm:AMGR} along the map $\boldsymbol{\Lambda}_{\bZ/\ell} \to \Lambda$.  \autoref{prop:compare-fh} implies that the $S^1$-actions agree, and \autoref{lem:finite-action-1} and \ref{lem:finite-action-2} imply that the maps $\boldsymbol{\Delta}_{\bZ/\ell} \to \Lambda$ and $\boldsymbol{\Lambda}_{\bZ/\ell} \to \Lambda$ agree after taking groupoid completions and are both the natural map $BC_\ell \to BS^1$. Hence the induced $C_\ell$-actions must also agree.
\end{proof}

Finally, using the cyclic bar construction over spectra, one can define the cyclotomic structure on the topological Hochschild homology. This shows that the cyclic traces over spectra can be lifted to cyclotomic spectra, by Nikolaus--Scholze \cite{Nikolaus-Scholze-Ontopologicalcyclichomology} and Ayala--Mazel-Gee--Rozenblyum \cite{AMGR-cyclotomic,AMGR-naive}:

\begin{Thm}[{\cite[Theorem B]{AMGR-cyclotomic}}]\label{thm:AMGR-cylotomic}
    Let $X$ be a small stable $\infty$-category (over $\Sp$). Then for $\ell \in \bN^\times$ there is a generalized Tate diagonal construction\footnote{Here we use the same convention as \cite{AMGR-cyclotomic,AMGR-naive}, where $\bN^\times$ means the monoid whose objects are positive natural numbers $\bN$ and monoidal structure is multiplication $\times$.}
    \[\varphi_\ell: \Tr(\Id_X) \to \Tr(\Id_X)^{\tau C_\ell}\]
    which assembles into a left-lax right action on $\Tr(\Id_X)$ and induces a canonical lift of $\Tr(\Id_X)$ to a cyclotomic spectrum.
\end{Thm}
\begin{RMK}
    When $\ell$ is a prime, the generalized Tate construction is equivalent to the classical Tate construction
    \[\Tr(\Id_X)^{\tau C_\ell} = \Tr(\Id_X)^{t C_\ell} := \mathrm{cofib}(\Tr(\Id_X)_{h C_\ell} \to \Tr(\Id_X)^{h C_\ell}).\]
    For an arbitrary $\ell \in \bN^\times$, the generalized Tate construction is the lower composite of the left Kan extension \cite[Remark 2.16]{AMGR-naive}
    \[\begin{tikzcd}
    \Sp^{BC_\ell} \ar[r,"(-)^{hC_\ell}"] \ar[d] & \Sp \\
    \Sp^{BC_\ell}/\{\Sigma^\infty+{C_\ell/H} \}_{H < C_\ell} \ar[ur]
    \end{tikzcd}\]
\end{RMK}
\begin{RMK}
   One may notice that the definition of the Tate diagonal package in \cite[Theorem A]{AMGR-cyclotomic} or \cite[Theorem B]{AMGR-naive} differs from \cite[Definition II.1.2]{Nikolaus-Scholze-Ontopologicalcyclichomology} where the latter only requires data for every single prime $p \in \bN^\times$. The reason more data is needed is that \cite{AMGR-cyclotomic,AMGR-naive} do not have the eventually connective assumption comparing to \cite[Theorem II.6.3]{Nikolaus-Scholze-Ontopologicalcyclichomology} (see \cite[Remark 0.11]{AMGR-naive}). Unfortunately, nor do we know the eventually connective assumption for the categories we are interested in later in the paper.
\end{RMK}
\begin{RMK}
    The cyclotomic structure is defined only on Hochschild homology over $\Sp$; see for example \cite[Theorem III.1.10]{Nikolaus-Scholze-Ontopologicalcyclichomology} \cite[Remark 3.9]{AMGR-naive}. Thus, it also seems difficult to write the Tate diagonal map formally in terms of the duality data.
\end{RMK}

\section{Category of sheaves and cyclic trace}\label{sec:cat-of-sheaf}

In this section, we will apply the results in \autoref{sec:trace} to study the Hochschild homologies of categories of sheaves, in both the topological and the algebro-geometric settings. We will write down the cyclic homologies without using explicit models of dg algebras.

\subsection{Admissible subcategories of presentable categories}

Let $\cE$ be a presentable locally rigid symmetric monoidal category (the assumption on local rigidity is because we need some discussion on adjunctions). Recall that for a $\cE$-linear category of $\cB$, we say a full $\cE$-linear subcategory is $\cA\hookrightarrow \cC$ is admissible if it is closed under both limit and colimit (no matter $\cE$-linear or not since $\cE$ is locally rigid). Let $\cB = \cC/\cA$. By the reflection theorem \cite{Reflection-Theorem-infty}, we have that $\cA$ is presentable $\cE$-linear and
\begin{equation}\label{eq:split-verider}
\cA \xrightarrow{\iota_*} \cC   \xrightarrow{j^!} \cB
\end{equation}
is a split Verdier sequence in $\Mod_\cE(\PrLst)$.
The fully faithful embedding $\iota_*: \cA \hookrightarrow \cC$ admits a left adjoint $\iota^*: \cC \to \cA$, and the quotient functor $j^!: \cC \twoheadrightarrow \cB$ also admits a fully faithful left adjoint $j_!: \cB \hookrightarrow \cC$ \cite[Lemma A.2.5]{Hermitian-K-theory}. Hence we have the following fiber sequence of functors
\begin{equation}L_\cB \longrightarrow \id\longrightarrow L_\cA,
\end{equation}
where $L_\cB = j_!j^!$ and $L_\cA = \iota_* \iota^*$. When $\cC$ is dualizable in $\Mod_\cE(\PrLst)$, we can often use \autoref{lemma: formula for dualizaing data} to deduce duality data for $\cA$ and $\cB$, by considering the idempotent $L_\cA$ (resp.~$L_\cB$) given by the composition
\[\cC \xrightarrow{\iota^*} \cA \xrightarrow{\iota_*} \cC, \quad (\text{resp.} \;\; \cC \xrightarrow{j^!} \cB \xrightarrow{j_!} \cC).
\]

\begin{Lemma}\label{lem: duality admissible category}
Let $\cA \hookrightarrow \cC \rightarrow \cB$ and $\cD \hookrightarrow \cC^\vee \rightarrow \mathcal{F}$ be split Verdier sequences in $\Mod_\cE(\PrLst)$. Suppose $(\cC, \cC^\vee, \epsilon_\cC, \eta_\cC)$ is a duality pair and $L_\cA^\vee = L_\cD$. Then the pair
\[\eta_\cA = (\iota_{\cD}^* \otimes \iota_{\cA}^*) \circ \eta_\cC, \quad \epsilon_\cA = \epsilon_\cC \circ (\iota_{\cA*} \otimes \iota_{\cD*})\]
exhibits $\cD$ as the dual of $\cA$, and the pair
\[
\eta_\cB = (j_{\mathcal{F}}^! \otimes j_{\cB}^!) \circ \eta_\cC, \quad \epsilon_\cB = \epsilon_\cC \circ (j_{\cB!} \otimes j_{\mathcal{F}!})
\]
exhibits $\mathcal{F}$ as the dual of $\cB$. 

Furthermore, under the identification $\Fun_\cE^L(\cA, \cA) = \cA^\vee \otimes \cA$, we have
\[L_\cA = \iota_\cD^*\iota_{\cD*} \otimes \iota_\cA^*\iota_{\cA*}, \quad L_\cB = j_{\mathcal{F}!}j_{\mathcal{F}}^! \otimes j_{\cB!}j_\cB^!.\]
\end{Lemma}

In the following subsections, we will see some examples from sheaf theory.

\subsection{Sheaves on manifolds with microsupport}
Let $X$ be a locally compact Hausdorff space, and we take a symmetric monoidal stable category $\cE\in \CAlg(\PrLst)$. We consider $\Sh(X;\cE)$, the $\cE$-linear category of $\cE$-valued sheaves. Then we have the $\cE$-linear 6-functor formalism on the category of sheaves \cite{6functor-infinity}. 

Let $\pi: X \to pt$ be the projection. We will write $1_X = \pi^*1_\cE$ for the constant sheaf on $X$. Let $i_Z: Z \hookrightarrow X$ be an embedding. We will write $1_Z = i_{Z!}1_Z$ for the standard sheaf on the subset $Z \subset X$. When $Z = \{x\}$, we will abuse notations and write $1_x = 1_{\{x\}}$.

Moreover, we also have the K\"unneth formula for the category of sheaves, cf.~\cite[Proposition~2.30]{6functor-infinity}. The following proposition is well-known, and a reference is {\cite[Proposition 3.1 \& 3.2]{Hochschild-Kuo-Shende-Zhang}}.

\begin{Prop}\label{prop: Künneth formula of sheaves}
The category $\Sh(X;\cE)$ is self-dual in $\Mod_{\cE}(\PrLst)$ with unit $\eta_X = \Delta_!\pi^*$ and counit $\epsilon_X = \pi_!\Delta^*$. Moreover, we have the K\"unneth formula
\[\Sh(X; \cE) \otimes_\cE \Sh(Y; \cE) \simeq \Sh(X\times Y;\cE).\]
Duality and K\"unneth formula together imply the equivalence
\[\Fun^\mathrm{L}_\cE(\Sh(X;\cE),\Sh(Y;\cE)) \simeq \Sh(X\times Y;\cE).\]
When $X = Y$, the identity functor $\Id_X$ is sent to $\eta_X = \Delta_!\pi^*1_\cE = \Delta_!1_{X}$.
\end{Prop}

Next, we assume $M$ is a smooth manifold and $\cE$ is compactly generated and rigid. 
To any object $F \in \Sh(M; \cE)$, one can associate a conic closed set $ss(F) \subset T^*M$ {\cite[Definition 5.1.2]{KS90}}, which satisfies the following triangular inequality together with other properties: for a fiber sequence $A \rightarrow B \rightarrow C$, we have $ss(A) \subset ss(B) \cup ss(C)$ \cite[Proposition~5.1.3]{KS90}.

\begin{RMK}Regarding microlocal sheaf theory in this setting, we remark that all arguments of \cite{KS90} work well provided we have the non-characteristic deformation lemma \cite[Proposition 2.7.2]{KS90} for all sheaves. When $\cE$ is compactly generated, the non-characteristic deformation lemma is proven for all hypersheaves \cite{Amicrolocallemma_infinitycat}, and then for all sheaves because of hypercompleteness for manifolds \cite[Thm.~7.2.3.6, Coro.~7.2.1.12]{HTT}.

The requirement of the coefficient category being compactly generated is a defect of the definition of microsupports given in \cite{KS90}. A possible definition for general coefficient categories was proposed by Efimov in \cite[Remark 4.24]{Efimov-K-theory} and was realized in \cite{Zhang_K_Remark}.    
\end{RMK}

We remark a further closure property of microsupports \cite[Exercise V.7]{KS90}, whose proof can be found in \cite[Proposition~3.4]{guillermouviterbo_gammasupport}. As a consequence, we have:

\begin{Prop}[{\cite[Proposition~6.26]{guillermouviterbo_gammasupport}}]\label{cor:sheaf-bicomplete}
For a projective system of sheaves $F_n\in \Sh(M;\cE)$ indexed by $n \in \bN$, we have 
\[ss\big(\varprojlim\nolimits_{n\in \bN}F_n\big) \subset \bigcap\nolimits_{N \geq 1}\overline{\bigcup\nolimits_{n\geq N} ss(F_n)}.\]
For an injective system of sheaves $F_n\in \Sh(M;\cE)$ indexed by $n \in \bN$, we have 
\[ss\big(\varinjlim\nolimits_{n\in \bN}F_n\big) \subset \bigcap\nolimits_{N \geq 1}\overline{\bigcup\nolimits_{n\geq N} ss(F_n)}.\]
\end{Prop}

\begin{Coro}\label{cor:split-verider-topology-sheaf}
Let $M$ be a manifold, $Z \subset T^*M$ be a closed conic subset and $U = T^*M \setminus Z$. Let $\Sh_Z(M; \cE)$ be the full subcategory of sheaves such that $ss(F) \subset Z$. Then 
\[\Sh_Z(M; \cE) \hookrightarrow \Sh(M; \cE) \rightarrow \Sh(M, U; \cE) := \Sh(M; \cE) / \Sh_Z(M; \cE)\]
is a split Verdier sequence in $\Mod_\cE(\PrLst)$.
\end{Coro}

By \autoref{prop: Künneth formula of sheaves} and \autoref{cor:split-verider-topology-sheaf}, we have a fiber sequence in $\Sh(X\times X;\cE)$
\[K_U\rightarrow 1_{\Delta_{X}}\rightarrow K_Z\]
representing the sequence of functors $L_{\Sh(M, U)}\rightarrow \id\rightarrow L_{\Sh_Z(M)}$ by convolution functors. The following lemma immediately follows from \autoref{lem: duality admissible category}.

\begin{Lemma}\label{lem:kernel as wrapping}
Let $\iota_{Z*}: \Sh_Z(M; \cE) \rightleftharpoons \Sh(M; \cE): \iota_Z^*$ be the inclusion and its left adjoint. Then we have
\[K_{Z} = \iota_{-Z \times Z}^*1_{\Delta_M}, \quad K_{U} = \fib(1_{\Delta_M} \to \iota_{-Z \times Z}^*1_{\Delta_M}).\]
Moreover, $K_Z$ is the unit that exhibits the duality $\Sh_Z(M; \cE)^\vee = \Sh_{-Z}(M; \cE)$, and $K_U$ is the unit that exhibits the duality $\Sh(M, U; \cE)^\vee = \Sh(M, -U; \cE)$.
\end{Lemma}

Therefore, applying \autoref{lemma: formula for dualizaing data} to the idempotent $L_{\Sh_Z(M)}$ or respectively $L_{\Sh(M, U)}$, we can write down the cyclic trace of the admissible subcategories explicitly using kernels, cf.~\cite[Section~3]{Hochschild-Kuo-Shende-Zhang} and \cite[Section~2]{zhang-cutoff}.

\begin{Coro}\label{cor:trace-sheaf}
Let $M$ be a manifold, $Z \subset T^*M$ be a closed conic subset and $U = T^*M \setminus Z$. Then the cyclic trace for $\Sh_Z(M; \cE)$ is given by
\[\Tr(\Id_{\Sh_Z(M)}, n)=\pi_{n!}\widetilde{\Delta}_n^* (K_{Z}\boxtimes \cdots\boxtimes K_{Z}), \]
where $\pi_{n} : X^n  \rightarrow \pt$ are projections, and $\Tilde{\Delta}_n$ are the twisted diagonals:
\begin{align*}
    \Tilde{\Delta}_n&: X^n \rightarrow X^{2n}, \qquad (x_1,\cdots,x_n)\mapsto (x_n,x_1,x_1,\cdots,x_{n-1},x_{n-1},x_n).
\end{align*}
The degeneracy maps are given by the morphism of sheaves $1_{\Delta_{M}} \to K_Z$, and the face maps are given by inverse of degeneracy maps. The cyclic trace for $\Sh(M, U; \cE)$ is given by
\[\Tr(\Id_{\Sh(M, U)}, n)=\pi_{n!}\widetilde{\Delta}_n^* (K_{U}\boxtimes \cdots\boxtimes K_{U} ). \]
The face maps are given by the morphism of sheaves $K_U\rightarrow 1_{\Delta_{X}}$, and the degeneracy maps are given by inverse of face maps as we explained in \autoref{rem:cyclic-object-map}.
\end{Coro}

By applying the results \autoref{prop:Gysin-finite-group} and \autoref{prop:Künneth-cyclic}, we have:

\begin{Coro}\label{cor: Gysin sheaf}
    Let $Z \subset T^*M$ be a closed conic subset. When $\cE$ is complex oriented, there is a Gysin sequence for the homotopy $S^1$ and $C_\ell$-quotients:
    \[\HH_*(\Sh_Z(M))_{hC_\ell} \to \HH_*(\Sh_Z(M))_{hS^1} \xrightarrow{[\ell](u)} \HH_*(\Sh_Z(M))_{hS^1}[2]\]
    and a Gysin sequence for the homotopy $S^1$ and $C_\ell$-fixed points:
    \[\HH_*(\Sh_Z(M))^{hS^1}[-2] \xrightarrow{[\ell](u)} \HH_*(\Sh_Z(M))^{hS^1} \to \HH_*(\Sh_Z(M))^{hC_\ell}\]
    where $u$ is the first Chern class in $H^2(BS^1; E)$ and $[\ell](u)$ is the $\ell$-series of the formal group law.
\end{Coro}

\begin{Coro}\label{cor: Künneth sheaf}
    Let $Z \subset T^*M$, $Y \subset T^*N$ be closed conic subsets. When $\cE$ is complex oriented, there is a K\"unneth isomorphism for the homotopy $S^1$-quotients:
    \[\HH_*(\Sh_{Z\times Y}(M \times N))_{hS^1} \simeq \HH_*(\Sh_Z(M))_{hS^1} \otimes _{E^{hS^1}} \HH_*(\Sh_Y(N))_{hS^1}[1],\]
    and a K\"unneth isomorphism for the homotopy $S^1$-fixed points within derived complete $\cE$-modules:
    \[\HH_*(\Sh_{Z\times Y}(M \times N))^{hS^1} \simeq \HH_*(\Sh_Z(M))^{hS^1} \widehat\otimes _{E^{hS^1}} \HH_*(\Sh_Y(N))^{hS^1}.\]
\end{Coro}

Finally, in the case when $\cE = \Sp$, we also define a cyclotomic structure on the trace of sheaf categories by \autoref{thm:AMGR-cylotomic}:
\begin{Coro}
Let $M$ be a manifold, $Z \subset T^*M$ be a closed conic subset and $U = T^*M \setminus Z$. For a $\ell \in \bN_{+}$, there exist natural morphisms
\[\varphi_\ell: \HH_*(\Sh_Z(M)) \to \HH_*(\Sh_Z(M))^{\tau C_\ell}\]
fitting into a left-lax commutative diagram, which induces a cyclotomic structure on $\HH_*(\Sh_Z(M)) = \Tr(\Id_{\Sh_Z(M)})$. Similarly, there exist morphisms
\[\varphi_\ell: \HH_*(\Sh(M, U)) \to \HH_*(\Sh(M, U))^{\tau C_\ell}\]
fitting into a left-lax commutative diagram, which induces a cyclotomic structure on $\HH_*(\Sh(M, U)) = \Tr(\Id_{\Sh(M, U)})$. 
\end{Coro}

We can also recover the $p$-th power construction and the $p$-fold cap product by \autoref{prop:power construction}:

\begin{Coro}Let $p$ be a prime number. Let $M$ be a manifold, $Z \subset T^*M$ be a closed conic subset. 
\begin{enumerate}
    \item[(1)] Let $E = D(\bZ/p\bZ)$. Then there is a $p$-fold cap operation
    \[\HH^*(\Sh_Z(M)) \to \Hom(\HH_*(\Sh_Z(M))^{tC_p}, \HH_*(\Sh_Z(M))^{tC_p}).\]
    \item[(2)] Let $\cE = \Sp$. Then there is also an $p$-fold cap operation
    \[\HH^*(\Sh_Z(M)) \to \Hom(\HH_*(\Sh_Z(M))^{tC_p}, \HH_*(\Sh_Z(M))^{tC_p}).\]
\end{enumerate}
\end{Coro}

\begin{eg}
We consider the case of $\Sh(M; \cE)$. Then by \cite[Proposition 3.2 \& Corollary 3.3]{Hochschild-Kuo-Shende-Zhang} or \cite[Theorem 0.2]{Efimov-K-theory}, we know that
\[\HH_*(\Sh(M; \cE)) = H_c^*(M; E), \quad \HH^*(\Sh(M; \cE)) = H^*(M, E).\]
In fact, by \cite[Theorem 0.2]{Efimov-K-theory}, we know that the circle action on $\HH_*(\Sh(M; \cE))$ is trivial; when $\cE = \Sp$, the cyclotomic structure on $\HH_*(\Sh(M; \Sp))$ is also trivial.

We can also consider the $p$-th power construction on the Hochschild cohomology. This gives a map
\[H^*(M; E) \to (H^*(M; E)^{\otimes p})^{tC_p}.\]
Note that the action of $(H^*(M; E)^{\otimes p})^{tC_p}$ reduces to an action of $H^*(M; E)^{tC_p}$ by $p$-fold multiplication. Now, its further composition, namely the following map
\[H^*(M; E) \to (H^*(M; E)^{\otimes p})^{tC_p} \to H^*(M; E)^{tC_p}\]
is the Tate-valued Frobenius map \cite[Definition IV.1.1]{Nikolaus-Scholze-Ontopologicalcyclichomology}. This composition admits a further factorization \cite[Corollary IV.1.4]{Nikolaus-Scholze-Ontopologicalcyclichomology}
\[H^*(M; E) \to (H^*(M; E)^{tC_p})^{h\operatorname{Aut}(C_p)} \to H^*(M; E)^{tC_p}.\]
When  $E = \bZ/p\bZ$, the first map gives exactly the Steenrod operation \cite[Theorem IV.1.15]{Nikolaus-Scholze-Ontopologicalcyclichomology}. When $E = KU$, the first map is the stable Adams operation \cite[Proposition IV.1.12]{Nikolaus-Scholze-Ontopologicalcyclichomology}.
\end{eg}

\begin{eg}Let $M$ be a closed manifold of dimension $d$. We consider the case of $\Loc(M; \cE)$. Then by Goodwillie--Jones \cite{Jones1987,Goodwillie85} and \cite[Section 5.1]{Brav-Dyckerhoff1} we know that
\[\HH_*(\Loc(M; \cE)) = H_*(LM; E), \quad \HH^*(\Loc(M; \cE)) = H_*(LM^{-TM}; E).\]
Moreover, the product on $\HH^*(\Loc(M; \cE))$ is given by the Chas--Sullivan loop product \cite{Cohen_2002} when $M$ is simply connected and $E$ is a discrete ring \cite{felixThomas2009string}. 

On the one hand, the circle action on $\HH_*(\Loc(M; \cE))$ is given by the loop rotation on $LM$; the cyclotomic structure on $\HH_*(\Loc(M; \cE))$ is given by the $p$-fold covering map of free loops $LM \to LM^{hC_p}$ by \cite[Proposition IV.3.2]{Nikolaus-Scholze-Ontopologicalcyclichomology}.

On the other hand, consider the $p$-th power construction on the Hochschild cohomology. This gives a map
\[H_*(LM^{-TM}, E) \to (H_*(LM^{-TM}, E)^{\otimes p})^{tC_p} \to \End(H_*(LM, E)^{tC_p}).\]
Assume that $M$ is $E$-orientable. Then we have a Calabi--Yau structure given by the cap product with a fundamental class of $M$ \cite{Brav-Rozenblyum}
\[\HH^*(\Loc(M)) \rightarrow \HH_*(\Loc(M))[-d].\]
We suspect that, under the isomorphism $H_*(LM^{-TM}; E) \simeq H_*(LM; E)[-d]$, the $p$-th power operation defined by the product on Hochschild cohomology
\[H_*(LM^{-TM}; E) \to (H_*(LM^{-TM}; E)^{\otimes p})^{tC_p} \to H_*(LM^{-TM}; E)^{tC_p}\]
is also induced by the $p$-fold covering map of free loops $LM \to LM^{hC_p}$.
\end{eg}

\subsection{Coherent and quasi-coherent sheaves}
We explain similar constructions for the categories of (quasi)coherent sheaves in algebraic geometry. Let $X$ be a derived stack satisfying appropriate assumptions \cite{Ben-Zvi-Francis-Nadler,Drinfeld-Gaitsgory,Gaitsgory-Rozenblyum,Six-FunctorScholze,stefanich2025dualizability_deriveds_stack}. We will consider the stable categories of quasi-coherent sheaves $\QCoh(X)$ and ind-coherent sheaves $\IndCoh(X)$ on the stack $X$\footnote{Under appropriate assumptions on $X$, we know $\QCoh(X)$ is the ind-completion of $\Perf(X)$, and $\IndCoh(X)$ is the ind-completion of $\Coh(X)$. However, neither result is true for general derived stacks \cite[Section 1.1]{Ben-Zvi-Francis-Nadler} \cite[Theorem 0.3.4]{Drinfeld-Gaitsgory} \cite[Theorem 1.0.4]{stefanich2025dualizability_deriveds_stack}.}. When $X$ is a regular perfect stack over a field $\Bbbk$, we have that $\Coh(X) = \Perf(X)$ and $\QCoh(X) = \IndCoh(X)$.

For quasi-coherent sheaves $\QCoh(X)$, we can prove duality and (consequently) K\"unneth formula in a number of settings. This was first proved by To\"en for quasi-compact separated schemes \cite[Section 8.3]{Toen_morita}, then Ben-Zvi--Francis--Nadler for perfect stacks \cite[Section 3]{Ben-Zvi-Francis-Nadler}, later Drinfeld--Gaitsgory for qca stacks that are eventually coconnective over fields of characteristic $0$ \cite[Section 4]{Drinfeld-Gaitsgory}, and more recently Stefanich for Noetherian stacks with affine stabilizers and separated diagonal that are well-stabilized \cite{stefanich2025dualizability_deriveds_stack}. 

In this setting, a partially defined 6-functor formalism is constructed in \cite[Section 3]{Ben-Zvi-Francis-Nadler}, \cite[Sections 1 \& 2]{Drinfeld-Gaitsgory} or \cite[Section 8.4]{Six-FunctorScholze} (for example, \cite[Proposition 3.10]{Ben-Zvi-Francis-Nadler} implies base change formula for the push-forward of perfect projections, and \cite[Corollary 1.4.5]{Drinfeld-Gaitsgory} implies base change formula for the push-forward of qca projections).

\begin{Prop}[{\cite[Corollary 1.8]{Toen_morita}, \cite[Theorem 1.2]{Ben-Zvi-Francis-Nadler}, \cite[Theorem 4.3.1]{Drinfeld-Gaitsgory}}]\label{prop:quasi-coh}
    Let $X$ and $Y$ be either perfect stacks over a ring $\Bbbk$, or qca algebraic stacks over a field $\Bbbk$ of characteristic $0$ that are eventually coconnective.
    The category $\QCoh(X)$ is self-dual in $\Mod_\Bbbk(\PrLst)$ with unit $\eta_X = \Delta_*\pi^*$ and counit $\epsilon_X = \pi_*\Delta^*$. Moreover, we have the K\"unneth formula
    \[\QCoh(X) \otimes_\Bbbk \QCoh(Y) \simeq \QCoh(X \times Y).\]
    Duality and K\"unneth formula together imply the equivalence
    \[\Fun^\mathrm{L}_\Bbbk(\QCoh(X), \QCoh(Y)) \simeq \QCoh(X \times Y).\]
    When $X = Y$, the identity functor $\Id_X$ is sent to $\eta_X = \Delta_*\pi^*\cO_{\Spec \Bbbk} = \Delta_*\cO_X$.
\end{Prop}

\begin{RMK}
    When $X$ is a perfect derived stack \cite{Ben-Zvi-Francis-Nadler}, $\QCoh(X)$ is also compactly generated. However, when $X$ is a qca derived stack \cite{Drinfeld-Gaitsgory} or well-stabilized Noetherian stack with affine stabilizers and separated diagonal \cite{stefanich2025dualizability_deriveds_stack}, $\QCoh(X)$ is dualizable but not known to be compactly generated. For a general laft derived stack \cite{Gaitsgory-Rozenblyum}, it is not true that $\QCoh(X)$ is dualizable.
\end{RMK}

For ind-coherent sheaves $\IndCoh(X)$, we also have duality and K\"unneth formula, but they are given by different functors, proved by Preygel for almost finitely presented derived stacks that are Noetherian, perfect, Deligne--Mumford with affine diagonal over perfect fields \cite[Appendix B]{Preygel}, Drinfeld--Gaitsgory for qca derived stacks \cite[Section 3 \& 4]{Drinfeld-Gaitsgory}, more generally Gaitsgory--Rozenblyum for locally almost finite type (laft) derived prestacks \cite[Chapter 5 Theorem 4.2.5]{Gaitsgory-Rozenblyum} \cite[Chapter 3 \S 5, ~\S 6]{Gaitsgory-RozenblyumII}, and Stefanich for well-stabilized Noetherian stacks with affine stabilizers and affine diagonals \cite[Theorem 1.0.4]{stefanich2025dualizability_deriveds_stack}.

In this setting, a complete 6-functor formalism is obtained in \cite[Part~II Chapter 5 Theorem 4.1.2 \& Part~III Chapter 8 Theorem 6.1.5]{Gaitsgory-Rozenblyum} for locally almost finite type ind-inf-schemes or \cite[Section 8.7]{Six-FunctorScholze} for almost finite type derived schemes, and to a large extent generalized to ind-geometric stacks \cite{CautisWilliams}.

    \begin{Prop}[{\cite[Proposition~B.3.2]{Preygel}, \cite[Theorem 3.3.5 \& Corollary 4.2.3]{Drinfeld-Gaitsgory}, \cite[Part~II Ch.~5 Thm.~4.2.5]{Gaitsgory-Rozenblyum}, \cite[ Ch.~3, Theorem~6.2.2 and ~\S6.2.6]{Gaitsgory-RozenblyumII}}]
    Let $X$ and $Y$ be laft ind-inf-schemes over a field $\Bbbk$. The category $\IndCoh(X)$ is self-dual in $\Mod_\Bbbk(\PrLst)$ with unit $\eta_X = \Delta_*\pi^!$ and counit $\epsilon_X = \pi_*\Delta^!$. Moreover, we have the K\"unneth formula
    \[\IndCoh(X) \otimes_\Bbbk \IndCoh(Y) \simeq \IndCoh(X \times Y).\]
    Duality and K\"unneth formula together imply the equivalence
    \[\Fun^\mathrm{L}_\Bbbk(\IndCoh(X), \IndCoh(Y)) \simeq \IndCoh(X \times Y).\]
    When $X = Y$, the identity functor $\Id_X$ is sent to $\eta_X = \Delta_*\pi^!\cO_{\Spec \Bbbk} = \Delta_*\omega_X$.
\end{Prop}

\begin{RMK}
When $X$ and $Y$ are smooth and proper perfect stacks, the Fourier--Mukai transformation theorem actually works for the small categories, by the results of Bondal--Larsen--Lunts \cite{Bondal-Larsen-Lunts}, To\"en \cite[Corollary 1.8]{Toen_morita} and Ben-Zvi--Francis--Nadler \cite[Theorem 1.2]{Ben-Zvi-Francis-Nadler}:
\[\Fun^{ex}_\Bbbk(\Perf(X), \Perf(Y)) \simeq \Perf(X) \otimes_\Bbbk \Perf(Y).\]
\end{RMK}

Similar to the topological setting, we can describe the cyclic homologies of admissible subcategories of coherent sheaves using a dg-algebra free model.

\begin{Coro}\label{cor:quasicoherent}
Let $X$ be a perfect stack over a ring $\Bbbk$ or a qca algebraic stack over a field $\Bbbk$ of characteristic $0$ and $\cA \subset \QCoh(X)$ be an admissible subcategory with integral kernel $K_\cA \in \QCoh(X\times X)$. Then the cyclic trace of $\cA$ is given by
\[\Tr(\Id_\cA, n) = \pi_{n*}\widetilde{\Delta}_n^*(K_\cA \boxtimes \cdots \boxtimes K_\cA).\]
\end{Coro}

\begin{Coro}\label{cor:coherent}
Let $X$ be a laft ind-inf-scheme over a perfect field $\Bbbk$  of characteristic $0$ and $\cA \subset \IndCoh(X)$ be an admissible subcategory with integral kernel $K_\cA \in \IndCoh(X\times X)$. Then the cyclic trace of $\cA$ is given by
\[\Tr(\Id_\cA, n) = \pi_{n*}\widetilde{\Delta}_n^!(K_\cA \boxtimes \cdots \boxtimes K_\cA).\]
\end{Coro}

\begin{RMK}
When $X$ is regular, we know that $\QCoh(X) = \IndCoh(X)$, so the description also works for ind-coherent sheaves. When $X$ is a smooth and proper variety, our description of Hochschild homology agrees with Kuznetsov \cite[Theorem~4.5 \& Proposition~8.1]{KuznetsovHochschild}.
\end{RMK}

\begin{eg}
For a smooth scheme $X$ over a field of characteristic $0$, one can interpret right (resp.~left) $D$-modules on $X$ as ind-coherent (resp.~quasi-coherent) sheaves on the de Rham stack $X_\mathrm{dR}$ \cite{Gaitsgory-Rozenblyum-crystal}. Recall that the functor $X \to X_\mathrm{dR}$ is the right Kan extension of the inclusion of reduced classical affine schemes into derived affine schemes \cite[Section 1.1]{Gaitsgory-Rozenblyum-crystal}:
\[X_\mathrm{dR}(S) = X(S^{cl,red}) \coloneqq X(H^0(S)/\mathrm{nilp}(H^0(S))).\]
Gaitsgory--Rozenblyum proved the following
\cite[Section 5.5]{Gaitsgory-Rozenblyum-crystal}:
Let $X$ be a classical scheme of finite type. Then there is a (canonical) equivalence
\[\IndCoh(X_\mathrm{dR}) \simeq D\mathrm{-mod}^r(X), \quad \QCoh(X_\mathrm{dR})\simeq D\mathrm{-mod}^l(X).\]
Since left and right $D$-modules are equivalent (for any laft pre-stacks) \cite[Proposition 2.2.4]{Gaitsgory-Rozenblyum-crystal}, we will no longer distinguish them and simply write $D\mathrm{-mod}(X)$.

The K\"unneth formula and (partially defined) 6-functor formalism for quasi-coherent sheaves then imply the relevant results for $D$-modules as well \cite{Gaitsgory-Rozenblyum-crystal} \cite[Corollary 8.3.4 \& 8.4.3]{Drinfeld-Gaitsgory}:
Let $X$ and $Y$ be qca algebraic stacks over a field $\Bbbk$ of characteristic $0$, and assume that $X$ is safe in the sense of \cite[Definition 10.2.2]{Drinfeld-Gaitsgory}. The category $D\mathrm{-mod}(X)$ is self-dual in $\Mod_\Bbbk(\PrLst)$ with unit $\Delta_{*,\mathrm{dR}}\pi^!$ and counit $\pi_{*,\mathrm{dR}}\Delta^!$. Moreover, we have the K\"unneth formula
\[D\mathrm{-mod}(X) \otimes_\Bbbk D\mathrm{-mod}(Y) \simeq D\mathrm{-mod}(X \times Y).\]
Duality and K\"unneth formula together imply the equivalence
\[\Fun_\Bbbk^\mathrm{L}(D\mathrm{-mod}(X), D\mathrm{-mod}(Y)) \simeq D\mathrm{-mod}(X \times Y).\]
Therefore, one can also consider admissible subcategories of $D$-modules and obtain a similar description on their Hochschild homology.
\end{eg}

\begin{RMK}
Hochschild invariants for coherent and quasi-coherent sheaves for schemes are closely related to the de Rham complex. For instance, if $X$ is smooth over a field of characteristic $0$, the well-known Hochschild--Konstant--Rosenberg isomorphism identifies the Hochschild homology $\HH_*(\QCoh(X))$ with differential forms on $X$ \cite{HKR,LodayCyclic,Yekutieli}, and the periodic cyclic homology $\HH_*(\QCoh(X))^{tS^1}$ with the periodic de Rham cohomology \cite{FeigenTsygan,LodayCyclic,toen_vezzosi_2011}. Over characteristic $p$, the Hochschild--Konstant--Rosenberg spectral sequence does not always degenerate \cite{AntieauBhattMathew}. Similar results for the matrix factorization category can also be obtained \cite{Preygel,Efimov-MF}.
When we take the topological Hochschild homology, then when $X$ is smooth over a perfect field of characteristic $p$, the associated graded of the BMS (Bhatt-Morrow-Scholze) filtration of $p$-complete periodic cyclic homology $\HH_*(\QCoh(X))^{tS^1}$ can be identified with the crystalline cohomology \cite{AntieauNikolaus}.

On the other hand, using derived algebraic geometry, To\"en--Vezzosi and Ben-Zvi--Nadler identified the Hochschild homology $\HH_*(\QCoh(X))$ with functions on the free loop space $LX$, and they identified the $S^1$-action on both sides \cite{toen_vezzosi_2011,toen_vezzosi_15,Ben-Zvi-Nadlerloopspaceandconnection}, which turned out to be equivalent to the de Rham cohomology on $X$. 

Here, what we do is to provide a dg-algebra free model for the Hochschild homology and cyclic homology of the coherent and quasi-coherent sheaves, in the spirit of similar constructions for non-equivariant Hochschild invariants \cite{KuznetsovHochschild,Ben-Zvi-Francis-Nadler,Preygel,Ben-Zvi-Nadler-trace}. (As explained by \cite{toen_vezzosi_2011,Ben-Zvi-Nadler-trace,Hoyois-Scherotzke-Sibilla}, one can define the equivariant lift by the cobordism hypothesis \cite{Lurie-cobordism}).
\end{RMK}

\section{Tamarkin category and cyclic trace}
In this section, we use the results of \autoref{sec:trace} to study the Hochschild homologies of the Tamarkin categories that arise from the categories of topological sheaves that were discussed in \autoref{sec:cat-of-sheaf}.
Using the results in \autoref{sec:gysin}, we will study $S^1$ and $C_\ell$-equivariant symplectic capacities.

\subsection{Tamarkin category of a point}
\label{Tamarkin point}

Let us write $\Sh_-(\bR)$ for the category of sheaves with nonpositive microsupport, and correspondingly $\Sh_+(\bR)$ for the sheaves with nonnegative microsupport.  In this section we study
\[\sT := \Sh(\bR) / \Sh_-(\bR).\]
Consider the left adjoint of the quotient $\sT \hookrightarrow \Sh_+(\bR)$. It realizes $\sT$ as the full subcategory on $\Sh_+(\bR)$ of objects with trivial compactly supported global sections \cite[Proposition 4.2]{Hochschild-Kuo-Shende-Zhang}. 

Let $s: \bR \times \bR \rightarrow \bR$ be the addition map.
Consider the map 
$\star: \Sh(\bR) \times \Sh(\bR) \rightarrow \Sh(\bR)$  which is given by the formula
\[ (G, F) \mapsto s_! (G \boxtimes F) = s_! (p_2^* G \otimes p_1^* F ).\]
This defines a symmetric monoidal structure on $\Sh(\bR)$, with unit $1_{0}$.  More generally, if we write the 
translation map on $\bR$ as $\nT_c(t)=t+c$, then 
one checks readily that $T_{c*} H \cong 1_{c} \star H$ for any $H \in \Sh(\bR)$. 

Since the symmetric monoidal product $\star$ preserves colimits, we find that $\Sh(\bR)$ is a closed symmetric monoidal category. We can describe the internal Hom $\sHom^\star$ by the formula
\[\sHom^\star(F,G)=p_{2*}\sHom(p_1^{*}F,s^!G),\quad F,G\in \Sh(\bR),\]
such that, for $F,G,H\in \Sh(\bR)$, 
\[\Hom_{\Sh(\bR)}(F\star G,H)=\Hom_{\Sh(\bR)}(F,\sHom^\star(G,H)).\]

\begin{Lemma}[{\cite[Corollary 4.14]{GS2014}}] \label{lem: convolution estimates}
The following estimates hold
\begin{align*}
    \star : &\,\Sh(\bR)\times \Sh_-(\bR) \rightarrow \Sh_-(\bR),\, \text{ or } \Sh_-(\bR)\times \Sh(\bR) \rightarrow \Sh_-(\bR);\\
   \sHom^\star : &\,\Sh(\bR)^{op}\times \Sh_-(\bR) \rightarrow \Sh_-(\bR),\, \text{ or } \Sh_-(\bR)^{op}\times \Sh(\bR) \rightarrow \Sh_-(\bR).
\end{align*}
\end{Lemma}

\autoref{lem: convolution estimates} implies that
the closed symmetric monoidal structure determined by $\star$
descends to $\sT$.  Note that in $\sT$, we have 
$[1_{\bR_{\geq 0}}]\cong [1_{0}]$, giving two useful expressions for
the unit.  
Let us note also that the inclusion $\sT \hookrightarrow \Sh(\bR), [F] \mapsto F \star 1_{\bR_{\geq 0}}$ is oplax monoidal with respect to $\star$.  Indeed: $ 1_{\bR_{\geq 0}}\star1_{\bR_{\geq 0}}=1_{\bR_{\geq 0}}$, so
\[([F]\star [G])\star 1_{\bR_{\geq 0}}=[F\star G] \star 1_{\bR_{\geq 0}}= F\star G \star 1_{\bR_{\geq 0}}=(F \star 1_{\bR_{\geq 0}}) \star (G\star 1_{\bR_{\geq 0}})\]
shows that $\sT \hookrightarrow \Sh(\bR)$ is a monoidal ideal.

\begin{Lemma}[{\cite[Lemma 4.4]{Hochschild-Kuo-Shende-Zhang}}]\label{lemma: hom computation}For $F,G\in \sT$, we have
\[\Hom_{\sT}(F,\nT_{c*}G)=\Hom_{\Sh(\bR)}(F\star 1_{\bR_{\geq 0}},\nT_{c*}G)=\Gamma_{[-c,\infty)}(\bR,\sHom^{\star}(F,G)).\]
\end{Lemma}

We mention that there is an equivalence $\Gamma_{[-,+\infty)}(\bR, -)$ from the Tamarkin category of point to torsion persistence modules \cite[Lemma 4.9]{Hochschild-Kuo-Shende-Zhang} or \cite[Proposition 4.22]{Efimov-K-theory}; in particular, for any $F \in \sT$, the system $\Gamma_{[c,+\infty)}(\bR; F)$, ${c\in\bR}$, forms an $\bR$-filtered complex that is a persistence module:

\begin{Lemma}[{\cite[Lemma 4.9]{Hochschild-Kuo-Shende-Zhang}, \cite[Proposition 4.22]{Efimov-K-theory}, \cite[Section 4.2]{APT-Tatsuki-Zhang}}]\label{lem:day convolution}
Let $\cE$ be a compactly generated stable category. There is a monoidal fully faithful embedding \[\sT \hookrightarrow \Fun(\larbR, \cE), \quad F \mapsto [c \mapsto \Gamma_{[c,+\infty)}(\bR, F)],\]
where the monoidal structure on $\sT$ is the sheaf convolution above, and the monoidal structure on $\Fun(\larbR, \cE)$ is given by the (derived) Day convolution 
\[F \star_{Day} G(c) = \varinjlim_{a + b \geq c} F(a) \otimes G(b).\]
\end{Lemma}

Finally, we record the observation that when $\cE$ is a compactly generated stable category with complex orientation in the sense of \autoref{def:complex orientation}, so is the Tamarkin category over $\cE$:

\begin{Lemma}\label{lem:complex orientation tamarkin}
Let $\cE$ be a compactly generated stable category with complex orientation. Then $\sT = \Sh(\bR; \cE)/\Sh_-(\bR; \cE)$ also has a complex orientation.
\end{Lemma}

\subsection{Tamarkin category of a domain}\label{sec:tamarkin domain}
Using the K\"unneth formula \autoref{prop: Künneth formula of sheaves}, we know that
\[\Sh(M; \sT) \simeq \Sh(M) \otimes \Sh(\bR)/\Sh_-(\bR)\\
\simeq \Sh(M \times \bR)/\Sh_{\{\tau\leq 0\}}(M \times \bR).\]
We can further define categories of sheaves valued in $\sT$ with microsupport conditions. They are called Tamarkin categories \cite{tamarkin2013,GS2014}.

\begin{Def}[Tamarkin category]
Let $U \subseteq T^*M$ and $Z = T^*M \setminus U$. Consider the maps $ T^*(M \times \bR) \xleftarrow{i} J^1 M \xrightarrow{\pi} T^*M$ and  define
\begin{align}
  \Sh_Z(M,\sT)\coloneqq \Sh_{\pi^{-1}(Z) \cup \{\tau \leq 0\}}(M \times \bR)/\Sh_{\{\tau \leq 0\}}(M \times \bR) = \{F\, |\, 
\pi(i^{-1} ss(F)) \subset Z\} /\Sh_{\{\tau \leq 0\}}(M \times \bR).
\end{align}
The notion of microsupport descends to define
a map $\underline{ss}$ from objects of $\sT(U)$ to closed coisotropic subsets of $U$, and we can write
\[\Sh_Z(M,\sT) = \{F\, |\, \underline{ss}(F) \subset Z\}.\]
We define the Tamarkin category for $U \subset T^*M$ to be
\begin{align}
\sT(U) \coloneqq \Sh(M; \sT) /\Sh_Z(M;\sT).
\end{align}
\end{Def}

Note that for
$j: V \hookrightarrow U$, there is an evident natural surjection 
$j^*: \sT(U) \to \sT(V)$
induced by the inclusion $i_*: \Sh_{T^*M \setminus U}(M; \sT) \hookrightarrow \Sh_{T^*M \setminus V}(M; \sT)$, which preserves both limits and colimits. 
Then, both $i_*$ and $j^*$ admit left and right adjoints $i^*(i^!)$ and $j_!(j_*)$ given by either the wrapping formula or, equivalently, the microlocal cut-off formula; see, for example \cite{Hochschild-Kuo-Shende-Zhang,zhang-cutoff}.

It is explained in \cite[Section 5]{Hochschild-Kuo-Shende-Zhang} that both $\cA=\Sh_Z(M,\sT)$ and $\cB=\sT(U)$ are $\sT$-linear admissible subcategories of $\Sh(M,\sT)=\sT(T^*M)$. As in loc.cit., we denote $Q_U=K_\cA$ and $P_U=K_\cB$.

\begin{Prop}[{\cite[Proposition~6.9]{Hochschild-Kuo-Shende-Zhang}}]\label{prop:unit-counit-tamarkin}
    Let $U \subset T^*M$ be an open subset. Then $\Sh_Z(M; \sT)$ is a dualizable $\sT$-linear category with
    \[\epsilon_{\Sh_Z(M; \sT)} = \pi_!^\sT \Delta_\sT^*, \; \eta_{\Sh_Z(M; \sT)} = Q_U = \iota_{(T^*M \setminus -U) \times (T^*M \setminus U)}^*1_{\Delta_M}^\sT.\]
    and $\sT(U)$ is a dualizable $\sT$-linear category with 
    \[\epsilon_{\sT(U)} = \pi_!^\sT \Delta_\sT^*, \; \eta_{\sT(U)} = P_U = \operatorname{fib}(1_{\Delta_M}^\sT \to \iota_{(T^*M \setminus -U) \times (T^*M \setminus U)}^*1_{\Delta_M}^\sT).\]
\end{Prop}

\subsection{Algebraic circle action on traces of Tamarkin categories}\label{sec:Tamarkin trace}
Therefore, the following results immediately follow from \autoref{thm:AF} and \autoref{prop:Gysin}.

\begin{Thm}\label{thm:S1-action}
    Let $U \subset T^*M$ be an open subset. Then the canonical functorial $S^1$-action on $\HH_*^\sT(\sT(U)) = \Tr^\sT(\Id_{\sT(U)})$ is defined by the cyclic object $\Tr(\Id_{\sT(U)}, -): \boldsymbol{\Lambda} \to \sT$ such that
    \begin{align}\label{eq:trace-circle}
    \Tr^\sT(\Id_{\sT(U)}, n) = \pi_{n!}^\sT\widetilde{\Delta}_{n,\sT}^*(P_U^{\boxtimes_\sT n}),
    \end{align}
    where the cyclic rotation is induced by the cyclic permutation $t_n: (M \times M)^{n} \to (M \times M)^{n}$ of order $n$, and the face maps $d_i$ and degeneracy maps $s_i$ are induced by the triangle identities of the unit and counit.
\end{Thm}
\begin{proof}
    The statement follows from \autoref{thm:AF} and \autoref{prop:unit-counit-tamarkin}.
\end{proof}
\begin{RMK}
    Using the functoriality of the Tamarkin category $\sT(U)$ with respect to inclusions and Hamiltonian homeomorphisms \cite{asanoike2022complete,guillermouviterbo_gammasupport}, we can immediately obtain the functoriality of the $S^1$-equivariant Hochschild homologies.
\end{RMK}

Similarly, considering the $C_\ell$-action on the Hochschild homology, the following results follow immediately from \autoref{prop:inclusion-in-circle} and  \autoref{prop:Gysin-finite-group}:

\begin{Thm}\label{thm:Z/ell-action}
    Let $U \subset T^*M$ be an open subset. The canonical functorial $C_\ell$-action on $\HH_*^\sT(\sT(U)) = \Tr^\sT(\Id_{\sT(U)})$ is defined by the cyclic object $\Tr(\Id_{\sT(U)}, -): \boldsymbol{\Delta}_{\bZ/\ell} \to \sT$ such that
    \begin{align}\label{eq:trace-finite-group}
    \Tr^\sT(\Id_{\sT(U)}, \ell) = \pi_{\ell!}^\sT\widetilde{\Delta}_{\ell\sT}^*(P_U^{\boxtimes_\sT \ell}),
    \end{align}
    where the cyclic rotation is induced by the cyclic permutation $t_\ell: (M \times M)^{\ell} \to (M \times M)^{\ell}$ of order $\ell$. 
\end{Thm}
\begin{proof}
    The statement follows from \autoref{prop:inclusion-in-circle} and \autoref{prop:unit-counit-tamarkin}.
\end{proof}
\begin{RMK}
    Using the functoriality of the Tamarkin category $\sT(U)$ with respect to inclusions and Hamiltonian homeomorphisms \cite{asanoike2022complete,guillermouviterbo_gammasupport}, we can immediately obtain the functoriality of the $C_\ell$-equivariant Hochschild homologies.
\end{RMK}

As a consequence, we also get natural $S^1$-actions on the linear dual of the Hochschild homology (namely the Chiu--Tamarkin complex \cite{Capacities2021,CyclicZHANG}).
We can easily identify the construction of $S^1$ and $C_\ell$-equivariant Chiu--Tamarkin complexes in \cite{zhang2020quantitative,CyclicZHANG} with the Hochschild cohomology of Tamarkin categories. We will use the notations in loc.~cit..

\begin{Prop}\label{prop:iso-with-chiu-tamarkin}
    Let $U \subset T^*M$ be an open subset and $\Bbbk$ be a discrete field. Then there is an isomorphism
    \[\Tr^\sT(\Id_{\sT(U)}, n) = \pi_{n!}\widetilde{\Delta}_{n}^*(s_{n!}P_U^{\boxtimes n}) =: F_n(U, \Bbbk),\]
    and for any $L \in \bR$, there is an isomorphism
    \begin{align*}
    i_L^*\Tr^\sT(\Id_{\sT(U)}, n) = i_L^*\pi_{n!}\widetilde{\Delta}_{n}^*(s_{n!}P_U^{\boxtimes n}) =: F_{L,n}(U, \Bbbk).
    \end{align*}
    Moreover, the $S^1$ and $C_\ell$-actions on both constructions are compatible.
    
    On the linear dual, there is also an isomorphism
    \begin{align*}
    \Tr^\sT(\Id_{\sT(U)}, n)^\vee = \sHom^\star(\pi_{n!}\widetilde{\Delta}_{n}^*(s_{n!}P_U^{\boxtimes n}), 1_{[0,+\infty)}) =: C_{n}(U, \Bbbk)[-\dim M].
    \end{align*}
    and for any $L \in \bR$, there is an isomorphism
    \begin{align*}
    \left(i_L^*\Tr^\sT(\Id_{\sT(U)}, n)\right)^\vee = \Hom(\pi_{n!}\widetilde{\Delta}_{n}^*(s_{n!}P_U^{\boxtimes n}), 1_{[L,+\infty)}) =: C_{L,n}(U, \Bbbk)[-\dim M].
    \end{align*}
    Moreover, the $S^1$ and $C_\ell$-actions on both constructions are compatible.
\end{Prop}
\begin{proof}
    The statement follows from comparing \autoref{eq:trace-circle} with \cite[Equation~(3.10)]{CyclicZHANG} and comparing \autoref{eq:trace-finite-group} with \cite[Equation~(2.1)]{CyclicZHANG} or \cite[Def.~2.13]{zhang2020quantitative}, plus the uniqueness of the microlocal kernels $P_U$ and $Q_U$. Note that in \cite[Section~3]{CyclicZHANG}, we chose a fibrant resolution of $Q_U$ of $c$-soft sheaves, but the mixed complexes are the same because of the equivalence between the stable $\infty$-category of (hyper)sheaves and the model category of dg-category of sheaves \cite{Cohn-dg}.
\end{proof}

We can deduce the Gysin sequence for cyclic homology of the Tamarkin category immediately by \autoref{prop:Gysin-finite-group}.

\begin{Thm}\label{thm: gysin Tamarkin}
    Let $U \subset T^*M$ be an open subset.    
    When $\cE$ is a complex oriented compactly generated stable category, there is a $C_\ell$-Gysin sequence between the homotopy $S^1$ and $C_\ell$ quotients (i.e.~cyclic homologies)
    \[\HH_*^\sT(\sT(U))_{hC_\ell} \to \HH_*^\sT(\sT(U))_{hS^1} \xrightarrow{[\ell](u)} \HH_*^\sT(\sT(U))_{hS^1}[2],\]
    and homotopy $S^1$ and $C_\ell$ fixed points (i.e.~negative cyclic homologies)
    \[\HH_*^\sT(\sT(U))^{hS^1}[-2] \xrightarrow{[\ell](u)} \HH_*^\sT(\sT(U))^{hS^1}  \to \HH_*^\sT(\sT(U))^{hC_\ell},\]
    where $u$ is the first Chern class in $H^2(BS^1; E)$ and $[\ell](u)$ is $\ell$-series of the formal group law in $E$.
\end{Thm}

We end the section by the K\"unneth formula for cyclic homology of the Tamarkin category which follows from \autoref{prop:Künneth-cyclic}.

\begin{Thm}\label{thm:Künneth tamarkin}
Let $U \subset T^*M$ and $W \subset T^*N$ be open subsets and $\cE$ be a complex oriented compactly generated stable category. Then there are fiber sequences
\begin{align*}\HH_*^\sT(\sT(U \times W))_{hS^1} &\xrightarrow{\Delta} \HH_*^\sT(\sT(U))_{hS^1} \star \HH_*^\sT(\sT(W))_{hS^1} \\
&\xrightarrow{u \otimes 1 - 1 \otimes u} \HH_*^\sT(\sT(U))_{hS^1} \star \HH_*^\sT(\sT(W))_{hS^1}[2],\end{align*}
and
\begin{align*}  \HH_*^\sT(\sT(U))^{hS^1} \widehat{\star} \HH_*^\sT(\sT(W))^{hS^1} [-2]&\xrightarrow{u \otimes 1 - 1 \otimes u} \HH_*^\sT(\sT(U))^{hS^1} \widehat{\star} \HH_*^\sT(\sT(W))^{hS^1} \\
&\xrightarrow{\Delta}\HH_*^\sT(\sT(U \times W))^{hS^1},\end{align*}
which implies that there is a K\"unneth isomorphism
\[\HH_*^\sT(\sT(U \times W))_{hS^1} \simeq \HH_*^\sT(\sT(U))_{hS^1} \star_{(1_\sT)^{hS^1}} \HH_*^\sT(\sT(W))_{hS^1}[1],\]
and
\[\HH_*^\sT(\sT(U \times W))^{hS^1} \simeq \HH_*^\sT(\sT(U))^{hS^1} \widehat{\star}_{(1_\sT)^{hS^1}} \HH_*^\sT(\sT(W))^{hS^1},\]

\end{Thm}
\begin{proof}
When $\cE$ has a complex orientation, we know that $\sT$ also has a complex orientation by \autoref{lem:complex orientation tamarkin}. Then the result follows from \autoref{prop:Künneth-cyclic}.
\end{proof}

\subsection{Additional structure on trace of the Tamarkin category}
Recall that we also have a duality between Hochschild homology and cohomology of the Tamarkin categories, proved in \cite[Proposition 7.2]{Hochschild-Kuo-Shende-Zhang}.

\begin{Prop}[{\cite[Proposition 7.2 \& Corollary 7.3]{Hochschild-Kuo-Shende-Zhang}}]\label{prop:CY}
Let $U \subset T^*M$ be an open subset and $\omega_M$ be the dualizing sheaf on $M$. Then we have
\[\HH^*_\sT(\sT(U), \omega_M^\sT) \simeq \HH^\sT_*(\sT(U))^\vee.\]
In fact, when $M$ is $E$-orientable, $\sT(U)$ is a proper Calabi--Yau category.
\end{Prop}

When $M$ is orientable, as a corollary of the cyclic Deligne conjecture \cite{Brav-Rozenblyum}, we have a framed $E_2$-structure on 
\[\HH^*_\sT(\sT(U)) \simeq \HH_*^\sT(\sT(U))^\vee[-\dim M],\]
which also induces an $L_\infty$-structure on $\HH^*_\sT(\sT(U))_{hS^1}[2]$ \cite{Westerland-Equiv-Operad}.

\begin{Prop}
Let $M$ be $E$-orientable and $U \subset T^*M$ be an open subset. Then $\HH^*_\sT(\sT(U))_{hS^1}[2]$ is a gravity algebra.

When $E$ is a discrete ring, on each connected component of the gravity operad, the operation is defined by
\begin{align*}
m_k: (\HH^*_\sT(\sT(U))_{hS^1})^{\otimes k} &\to (\HH^*_\sT(\sT(U)))^{\otimes k}[-k] \\
&\to \HH^*_\sT(\sT(U))[-k] \to \HH^*_\sT(\sT(U))_{hS^1}[-k].
\end{align*}
where the first and third maps come from the Gysin sequence, and the second map comes from the product structure.
\end{Prop}
\begin{proof}
The proof is completely analogous to \cite[Section 7.4]{Westerland-Equiv-Operad}. Consider the $E_2^\textit{fr}$-algebra structure:
\[H_*(E_2^\textit{fr}(k), E) \otimes \HH^*_\sT(\sT(U))^{\otimes k} \to \HH^*_\sT(\sT(U)).\]
We now unpack how this induces a gravity algebra structure on the $S^1$-quotient following \cite[Section 7.2]{Westerland-Equiv-Operad}. First, taking the homotopy quotient, we have that
\[(H_*(E_2(k), E) \otimes \HH^*_\sT(\sT(U))^{\otimes k})_{hS^1} \to \HH^*_\sT(\sT(U))_{hS^1}.\]
Note that $S^1$ acts on $E_2(k)$ and $\HH^*_\sT(\sT(U))$. We can consider the transfer map in the Gysin sequence with respect to the inclusion $S^1 \hookrightarrow (S^1)^{k+1}$, this induces
\begin{align*}
H_*(E_2(k), E)_{hS^1} \otimes (\HH^*_\sT(\sT(U))_{hS^1})^{\otimes k} &\xrightarrow{\sim} (H_*(E_2(k), E) \otimes \HH^*_\sT(\sT(U))^{\otimes k})_{h(S^1)^{k+1}} \\
&\to (H_*(E_2(k), E) \otimes \HH^*_\sT(\sT(U))^{\otimes k})_{hS^1}[-k].
\end{align*}
When we restrict to the degree zero part of $H_*(E_2(k)_{hS^1}, E)$, the above transfer map agrees with the transfer map in the Gysin sequence
\begin{align*}
(\HH^*_\sT(\sT(U))_{hS^1})^{\otimes k} \xrightarrow{\sim}(\HH^*_\sT(\sT(U))^{\otimes k})_{h(S^1)^k} \to \HH^*_\sT(\sT(U))^{\otimes k}[-k].
\end{align*}
Composing with the multiplication and the quotient map, this is exactly the formula in the proposition.
\end{proof}

\begin{RMK}By \cite[Corollary 2.4]{Westerland-Equiv-Operad}, there is a similar algebraic structure on $\HH^*_\sT(\sT(U))^{hS^1}$, which is also deeply related to gravity algebra by \cite[Corollary 2.8]{Westerland-Equiv-Operad}.
\end{RMK}

\begin{RMK}
Siegel \cite{SiegelHigherCapacities} defined higher symplectic capacities using the $L_\infty$-structure on the $S^1$-equivariant symplectic cohomology and the $L_\infty$-morphism $\mathrm{SH}_+^{S^1}(U) \to \Bbbk[u]$ which induces a dg coalgebra map on the bar construction (the Chevalley--Eilenberg cochain complex)
\[B\mathrm{SH}_+^{S^1}(U) \to B(\Bbbk[u]).\]
Over a field of characteristic $0$, a gravity algebra in particular induces operations of an $L_\infty$-algebra. Based on the above proposition, one can similarly define higher symplectic capacities by formally taking the bar construction (the Chevalley--Eilenberg construction).
\end{RMK}

\bibliographystyle{bingyu}
\clearpage
\phantomsection
\bibliography{bibtexarxiv}

\end{document}